\documentclass[10pt]{amsart}
\usepackage{amsmath,amssymb,amsthm,amscd,
graphics,axodraw,verbatim}
\usepackage[dvips]{color}

\usepackage[stdtext,vcentermath]{myyoungtab}

\newtheorem{theorem}{Theorem}
\newtheorem{lemma}[theorem]{Lemma}

\newtheorem{proposition}[theorem]{Proposition}

\numberwithin{theorem}{section}

\theoremstyle{definition}
\newtheorem{definition}[theorem]{Definition}
\newtheorem{example}[theorem]{Example}

\newtheorem*{acknowledgments}{Acknowledgments}

\numberwithin{equation}{section}

\newcommand{\ot}{\otimes}

\newcommand{\Z}{{\mathbb Z}}

\newcommand{\J}{{\mathcal J}}

\newcommand{\Pth}{{\mathcal P}}

\newcommand{\Aff}{{\rm Aff}}
\newcommand{\wt}{{\rm wt}}

\def\et#1{\tilde{e}_{#1}}
\def\ft#1{\tilde{f}_{#1}}

\def\veps{\varepsilon}
\def\vphi{\varphi}
\def\ot{\otimes}

\begin{document}

\begin{center}
{\Large{\bf Combinatorial Bethe ansatz  \\
 and \\
Generalized periodic box-ball system}}
\vspace{10mm}\\
{\large Atsuo Kuniba and Reiho Sakamoto}
\vspace{3mm}
\end{center}

\begin{quotation}
{\small
A}{\tiny BSTRACT:}
{\small
We reformulate the Kerov-Kirillov-Reshetikhin (KKR) map 
in the combinatorial Bethe ansatz 
from paths to rigged configurations by introducing 
local energy distribution in crystal base theory.
Combined with an earlier result on the inverse map,
it completes the crystal interpretation of 
the KKR bijection for $U_q(\widehat{\mathfrak{sl}}_2)$.
As an application, we solve 
an integrable cellular automaton, 
a higher spin generalization of the periodic box-ball system,  
by an inverse scattering method and obtain the solution of the 
initial value problem in terms of 
the ultradiscrete Riemann theta function.}
\end{quotation}

\section{Introduction}\label{sec:intro}

The Kerov-Kirillov-Reshetikhin (KKR) bijection \cite{KKR,KR} 
is a combinatorial version of the Bethe ansatz.
It gives a one to one correspondence between 
rigged configurations and highest paths, which are 
combinatorial analogues of the 
Bethe roots and the associated Bethe vectors 
in integrable spin chains\footnote{
The original KKR bijection concerns semistandard tableaux
rather than highest paths. The KKR bijection 
in this paper is to be understood 
as the composition of the original one 
with the Robinson-Schensted-Knuth correspondence between 
semistandard tableaux and highest paths.}.
The relevant problem of state counting 
stemmed from Bethe's original work \cite{Be},
was developed further in the KKR theory, 
and has been formulated 
as the $X=M$ conjecture for arbitrary affine Lie algebra \cite{HKOTT}. 
See \cite{KSS, Sch1} for a recent status.

In \cite{KOSTY, Sa}, 
the KKR map $\phi^{-1}$ from rigged configurations
to highest paths was identified with a certain composition of 
combinatorial $R$ in crystal base theory \cite{K,KMN,NY}. 
It provided a long sought representation theoretical 
meaning with $\phi^{-1}$ and opened a connection with the 
integrable cellular automata 
called the box-ball system \cite{TS,TM}
and its generalizations \cite{HKT,FOY,HHIKTT}.
They are identifed with solvable vertex models \cite{Ba} 
associated with the quantum group $U_q$ at $q=0$. 
In this context, the KKR theory is regarded as 
the inverse scattering formalism of the generalized 
box-ball systems, where the rigged configurations and 
$\phi^{-1}$ play the roles of 
scattering data and inverse scattering transform, 
respectively.
Precise descriptions are available either in 
Proposition 2.6 in \cite{KOSTY},  
section 3.2 and appendix E in \cite{KSY}, and 
Lemma \ref{s:lem:haiesutoka} in this paper.

In this paper we study two closely related problems
concerning $U_q(\widehat{\mathfrak{sl}}_2)$ case.
In the first part (section \ref{sec:led}), 
we give a crystal theoretical interpretation of 
the opposite KKR map $\phi$ from paths to rigged configurations.
It is done by introducing the {\em local energy distribution}
of paths, which provides a bird's-eye view 
of the whole combinatorial procedures involved in the KKR algorithm.
In terms of generalized box-ball systems,  
$\phi$ is a direct scattering map and separates the 
dynamical degrees of freedom into action-angle variables,
which are amplitudes and phase of solitons.
The local energy distribution makes it possible to grasp 
these data from a global viewpoint. See Example \ref{ex:big}. 
Together with the earlier result on $\phi^{-1}$,
we complete the crystal interpretation of the KKR bijection
$\phi^{\pm 1}$ for $U_q(\widehat{\mathfrak{sl}}_2)$.

The results mentioned so far are concerned with 
generalized box-ball systems on (semi) infinite lattice. 
In the second part of the paper (section \ref{sec:pbbs}), 
we launch the inverse scattering formalism in the periodic case.
This was achieved in \cite{KTT} 
for the simplest spin $\frac{1}{2}$ system 
called the periodic box-ball system \cite{YT,YYT}, 
and subsequently in \cite{MIT}.
Here we treat the general spin $\frac{s}{2}$ case
based on the crystal base theory.
Here is an example of time evolution $(T_4)$ of an $s=3$ case.
\begin{center}
$ t=0\,:\;122 \cdot 122 \cdot 112 \cdot 112 \cdot 111 \cdot 122 \cdot 111 \cdot 111 \cdot 112 $\\
 $t=1\,:\;112 \cdot 112 \cdot 122 \cdot 122 \cdot 112 \cdot 111 \cdot 122 \cdot 111 \cdot 111 $\\
 $t=2\,:\;111 \cdot 112 \cdot 112 \cdot 112 \cdot 122 \cdot 122 \cdot 111 \cdot 122 \cdot 111 $\\
 $t=3\,:\;111 \cdot 111 \cdot 112 \cdot 112 \cdot 112 \cdot 112 \cdot 222 \cdot 111 \cdot 122 $\\
 $t=4\,:\;122 \cdot 111 \cdot 111 \cdot 112 \cdot 112 \cdot 112 \cdot 111 \cdot 222 \cdot 112$ \\
 $t=5\,:\;112 \cdot 222 \cdot 111 \cdot 111 \cdot 112 \cdot 112 \cdot 112 \cdot 111 \cdot 122$ \\
 $t=6\,:\;122 \cdot 111 \cdot 222 \cdot 112 \cdot 111 \cdot 112 \cdot 112 \cdot 112 \cdot 111$ \\
 $t=7\,:\;111 \cdot 122 \cdot 111 \cdot 122 \cdot 122 \cdot 111 \cdot 112 \cdot 112 \cdot 112$ \\
 $t=8\,:\;112 \cdot 111 \cdot 122 \cdot 111 \cdot 112 \cdot 222 \cdot 111 \cdot 112 \cdot 112$ \\
 $t=9\,:\;112 \cdot 112 \cdot 111 \cdot 122 \cdot 111 \cdot 111 \cdot 222 \cdot 112 \cdot 112$ \\
\hspace{-0.1cm}$t=10\,:\;112 \cdot 112 \cdot 112 \cdot 111 \cdot 122 \cdot 111 \cdot 111 \cdot 122 \cdot 122$ \\
\hspace{-0.1cm}$t=11\,:\;122 \cdot 122 \cdot 112 \cdot 112 \cdot 111 \cdot 122 \cdot 111 \cdot 111 \cdot 112$\\
\end{center}
A local spin 
$\frac{s}{2}$ state is an $s$ array of $1$ and $2$
which are arranged not to decrease to the right.
Each local state is regarded as a capacity $s$ box.
Local states, say $111, 112, 122$ and $222$ for $s=3$, 
represent an empty box and those containing 
$1,2$ and $3$ balls, respectively. 
An array of such local states are called paths.
The above paths are of length 9.

A path of length $L$ can naturally be viewed as an
element of $B_s^{\otimes L}$, the tensor product of the 
crystal $B_s$ of the $s$-fold symmetric tensor 
representation of $U_q(\widehat{\mathfrak{sl}}_2)$. 
A wealth of notions and combinatorial operations on $B_s$ are
provided by the crystal base theory.
We make use of them to characterize a certain class of paths that are
invariant under extended affine Weyl group 
${\widehat W}(A^{(1)}_1)$ and 
the commuting family of invertible time evolutions
$\{T_l\}$.
This is an important 
non-trivial step characteristic to the $s>1$ situation.    
We introduce action-angle variables which correspond to 
those paths bijectively and linearize the dynamics.
These features are integrated in Theorem \ref{th:cd}. 
As corollaries of it, 
generic period and a counting formula of the paths are obtained 
in terms of conserved quantities in (\ref{eq:lcm}) and 
(\ref{eq:omega}), respectively.
For example (\ref{eq:lcm}) tells that 
the period of the above paths under $T_4$ is indeed $11$.
(Notice that the $t=0$ and $t=11$ paths are the same.)
These results agree with the conjecture 
in the most general setting \cite{KT}.
The initial value problem is solved either by 
a combinatorial algorithm or by an explicit formula 
(\ref{eq:sol}), (\ref{eq:xyt}) involving the 
ultradiscrete Riemann theta function (\ref{eq:urt}), 
generalizing the $s=1$ results in \cite{KS1,KS2}.
These expressions follow rather straightforwardly from 
the ultradiscrete tau function studied in \cite{KSY}.
For the background idea of ultra-discretization and 
relevant issues in tropical geometry, see \cite{TTMS} and \cite{MZ}.

Several characteristic features in 
quasi-periodic solutions to soliton equations \cite{DT,DMN}
will be demonstrated in the ultradiscrete setting. 
In particular our action-angle variables live in the set 
(\ref{eq:av}) which is an ultradiscrete analogue of 
the Jacobi variety. 
For a reduced case (\ref{eq:js}) with $s=1$, 
the underlying tropical hyperelliptic curve has been identified 
recently \cite{IT}. 
The action-angle variables are essentially solutions of the 
string center equation, which is a version of  
the Bethe equation at $q=0$ \cite{KN}.  
In this sense, the inverse scattering formalism in this paper 
connects the Bethe ans\"atze at $q=1$ \cite{KR,KKR} 
and $q=0$ \cite{KN} to 
the algebraic geometry techniques of soliton theory 
at a combinatorial level.

Our crystal interpretation of the KKR map $\phi$
has stemmed from an attempt to formulate the direct scattering 
map in the generalized periodic box-ball system.
In fact, we will show in section \ref{ssec:dst} that 
the idea of local energy distribution 
is efficient also in the periodic setting.

The paper is organized as follows.
In section \ref{sec:led}, the KKR map $\phi$ is
identified with a procedure based on the local energy distribution
in Theorem \ref{s:main}.
We illustrate it along a few instructive examples.
The proof will be given in \cite{Sa2}. 
Section \ref{sec:pbbs} is devoted to 
the generalized periodic box-ball system.
Section \ref{sec:sum} is a summary.
Appendix \ref{sec:cry} recalls the 
basic facts on crystal base theory \cite{K,KMN,NY}.
Appendix \ref{sec:kkr} is an exposition of the 
KKR bijection including the non-highest case \cite{DS,Sch1}.

\section{Local energy distribution and the KKR bijection}\label{sec:led}
In this section, we reformulate the combinatorial
procedure of the KKR map $\phi$ in terms of the energy functions
of crystal base theory.
See Appendix \ref{sec:cry} for the basic facts on 
crystal base theory.  
Consider the relation 
\begin{equation*}
a\otimes b_1\simeq b_1'\otimes a'
\end{equation*}
and the energy function
$e_1=H(a\otimes b_1)$ under the combinatorial $R$.
We depict them by the vertex diagram:
\begin{center}
\unitlength 13pt
\begin{picture}(4,4)
\put(0,2.0){\line(1,0){3.2}}
\put(1.6,1.0){\line(0,1){2}}
\put(-0.6,1.8){$a$}
\put(1.4,0){$b_1'$}
\put(1.4,3.2){$b_1$}
\put(3.4,1.8){$a'$}
\put(0.8,2.3){$e_1$}
\put(4.1,1.7){.}
\end{picture}
\end{center}
Successive applications of the combinatorial $R$ 
\begin{equation*}
a\otimes b_1\otimes b_2\simeq b_1'\otimes a'\otimes b_2
\simeq b_1'\otimes b_2'\otimes a'',
\end{equation*}
with  $e_2=H(a'\otimes b_2)$ is 
expressed by joining two vertices:
\begin{center}
\unitlength 13pt
\begin{picture}(8,4)
\multiput(0,0)(4.2,0){2}{
\put(0,2.0){\line(1,0){3.2}}
\put(1.6,1.0){\line(0,1){2}}
}
\put(-0.6,1.8){$a$}
\put(1.4,0){$b_1'$}
\put(1.4,3.2){$b_1$}
\put(3.4,1.8){$a'$}
\put(5.6,3.2){$b_2$}
\put(5.6,0){$b_2'$}
\put(7.6,1.8){$a''$}
\put(0.7,2.3){$e_1$}
\put(4.9,2.3){$e_2$}
\put(8.3,1.7){.}
\end{picture}
\end{center}

Given a path $b=b_1\otimes b_2\otimes\cdots\otimes b_L$,
its local energy ${\mathcal E}_{l,j}$ is defined by
${\mathcal E}_{l,j}:=H(u_l^{(j-1)}\otimes b_j)$,
where $u_l^{(j-1)}$ is specified by the following diagram
with the convention $u_l^{(0)}=u_l$ (\ref{eq:ul}).
\begin{center}
\unitlength 13pt
\begin{picture}(22,6.2)(0,-1)
\multiput(0,0)(5.8,0){2}{
\put(0,2.0){\line(1,0){4}}
\put(2,0){\line(0,1){4}}
}
\put(-0.9,1.8){$u_l$}
\put(0.5,2.3){${\mathcal E}_{l,1}$}
\put(1.7,4.2){$b_1$}
\put(1.7,-0.8){$b_1'$}
\put(4.2,1.8){$u_l^{(1)}$}
\put(6.3,2.3){${\mathcal E}_{l,2}$}
\put(7.5,4.2){$b_2$}
\put(7.5,-0.8){$b_2'$}
\put(10.0,1.8){$u_l^{(2)}$}
\multiput(11.5,1.8)(0.3,0){10}{$\cdot$}
\put(14.7,1.8){$u_l^{(L-1)}$}
\put(17,0){
\put(0,2.0){\line(1,0){4}}
\put(2,0){\line(0,1){4}}
}
\put(17.4,2.3){${\mathcal E}_{l,L}$}
\put(18.7,4.2){$b_L$}
\put(18.7,-0.8){$b_L'$}
\put(21.2,1.8){$u_l^{(L)}$}
\end{picture}
\end{center}
We set ${\mathcal E}_{0,j}=0$ for all $1\leq j\leq L$.
We define ${\mathcal T}_l$ and ${\mathcal E}_l$ by 
${\mathcal T}_l(b)=b_1'\otimes b_2'\otimes\cdots\otimes b_L'$ and 
\begin{equation}\label{eq:eli}
{\mathcal E}_l:=\sum_{j=1}^L {\mathcal E}_{l,j}.
\end{equation}
In other words, $u_l[0]\otimes b\stackrel{R}{\simeq}{\mathcal T}_l(b)
\otimes u_l^{(L)}[{\mathcal E}_l]$, where we have omitted
modes for $b$ and ${\mathcal T}_l(b)$.

Given a path $b=b_1\otimes b_2\otimes\cdots\otimes b_L$
($b_i\in B_{\lambda_i}$), 
we always have $u_l^{(L+\Lambda)}=u_l$ for any $l$ for 
a modified path 
$b'=b\otimes \fbox{1}^{\,\otimes\Lambda}$
if $\Lambda >\lambda_1+\cdots +\lambda_L$.
In such a circumstance, ${\mathcal E}_l({\mathcal T}_k(b'))={\mathcal E}_l(b')$ is known to hold
(Theorem 3.2 of \cite{FOY}, section 3.4 of \cite{HHIKTT}).
Namely the sum ${\mathcal E}_l$ is a conserved
quantity of the box-ball system on semi-infinite lattice.
Here we need more detailed information such as ${\mathcal E}_{l,j}$.

\begin{lemma}
For a path $b=b_1\otimes b_2\otimes\cdots\otimes b_L$,
we have ${\mathcal E}_{l,j}-{\mathcal E}_{l-1,j}=0$ or 1,
for all $l>0$ and for all $1\leq j\leq L$.
\end{lemma}
\begin{proof}
When $l=1$, this is clear from the definition ${\mathcal E}_{0,j}=0$
and the fact $H(x\otimes y)=0$ or 1 for any $x\in B_1$.

Now we investigate possible values for ${\mathcal E}_{l,j}-{\mathcal E}_{l-1,j}$.
We show that the difference between
tableaux for $u_{l}^{(j)}$ and $u_{l-1}^{(j)}$
is only one letter, namely, if $u_{l-1}^{(j)}=(x_1,x_2)$,
then $u_{l}^{(j)}=(x_1+1,x_2)$ or
$u_{l}^{(j)}=(x_1,x_2+1)$.
We show the claim by induction on $j$.
For $j=0$, it is true because
$u_{l-1}^{(0)}=u_{l-1}=(l-1,0)$ and $u_{l}^{(0)}=u_{l}=(l,0)$.
Suppose that the above claim holds for all $j<k$
for some $k$.
In order to compare $u_{l-1}^{(k)}$ and $u_l^{(k)}$,
consider the relations
$u_{l-1}^{(k-1)}\otimes b_{k}\simeq b_{l-1,k}'\otimes u_{l-1}^{(k)}$
and
$u_l^{(k-1)}\otimes b_{k}\simeq b_{l,k}'\otimes u_l^{(k)}$.
By the assumption, the difference between $u_{l-1}^{(k-1)}$
and $u_{l}^{(k-1)}$ is one letter.
Recall that in calculating the combinatorial $R$ by the graphical rule
(section \ref{subsec:r}),
order of making pairs is arbitrary.
Therefore, in $u_l^{(k-1)}\otimes b_{k}$,
first we can make all pairs that appear in
$u_{l-1}^{(k-1)}\otimes b_{k}$,
and then make the remaining one pair.
This means the difference of number of unwinding pairs,
i.e., ${\mathcal E}_{l,k}-{\mathcal E}_{l-1,k}$ is 0 or 1.
To make the induction proceed, note that this fact
means the difference between $u_{l-1}^{(k)}$
and $u_{l}^{(k)}$ is also one letter.
\end{proof}

Let $b=b_1\otimes\cdots\otimes b_L
\in B_{\lambda_1}\otimes\cdots\otimes
B_{\lambda_L}$ be an arbitrary (either highest or not) path.
Set $N={\mathcal E}_1(b)$.
We determine the pair of numbers
$(\mu_1,r_1)$, $(\mu_2,r_2)$, $\cdots$, $(\mu_N,r_N)$
by Step (i)--(iv).
\begin{enumerate}
\item
Draw a table containing $({\mathcal E}_{l,j}-{\mathcal E}_{l-1,j}=0,1)$
at the position $(l,j)$, i.e., at the $l$ th row and
the $j$ th column.
We call this table local energy distribution.
\item
Starting from the rightmost 1 in the $l=1$ st row,
pick one 1 from each successive row.
The one in the $(l+1)$ th row must be weakly right of
the one selected in the $l$ th row.
If there is no such 1 in the $(l+1)$ th row,
the position of the lastly picked 1 is called $(\mu_1,j_1)$.
Change all the selected 1 into 0.
\item
Repeat Step (ii) $(N-1)$ times to further determine
$(\mu_2,j_2)$, $\cdots$, $(\mu_N,j_N)$
thereby making all 1 into 0.
\item
Determine $r_1,\cdots,r_N$ by
\begin{equation}\label{s:eq:rigging}
r_k=\sum_{i=1}^{j_k-1}\min (\mu_k,\lambda_i)
+{\mathcal E}_{\mu_k,j_k}-2\sum_{i=1}^{j_k}{\mathcal E}_{\mu_k,i}.
\end{equation}
\end{enumerate}

One may replace the procedure (ii) by

\vspace{0.2cm}
(ii)'\; 
Starting from any one of the lowest $1$, pick one $1$ from 
each preceding row. The one \\
\phantom{(ii)'\; \;\;\;}
in the $(l-1)$th row must be weakly 
left and nearest of the one selected in the $l$th 
\phantom{(ii)'\; \;\;\;\,}row.
The position of the firstly picked $1$ is called $(\mu_1, j_1)$.
Change all the selected $1$ 
\phantom{(ii)'\; \;\;\;\,}into $0$. 

\vspace{0.2cm}
Our main result in this section is the following theorem, which 
gives a crystal theoretic
reformulation of the KKR map $\phi$.

\begin{theorem}\label{s:main}
The above procedure (i)--(iv) is well defined and 
$(\lambda ,(\mu ,r))$ coincides with
the (unrestricted) rigged configuration $\phi (b)$.
The procedure (i), (ii)',(ii),(iv) is also well defined and 
leads to the same 
rigged configuration up to a permutation of 
$(\mu_k, j_k)$'s.

\end{theorem}

The proof will be given in \cite{Sa2}.

\begin{example}\label{s:ex:mainth}
Consider the path which will also be treated in Example \ref{s:ex:kkr}:
\begin{equation*}
b=
\fbox{1111}\otimes\fbox{11}\otimes\fbox{22}\otimes
\fbox{12}\otimes\fbox{2}\otimes\fbox{122}\otimes
\fbox{122}\otimes\fbox{1112}
\end{equation*}
According to Step (i),
the local energy distribution is given in
the following table
($j$ stands for column coordinate of the table).
\begin{center}
\begin{tabular}{l|cccccccc}
                 &1111&11&22&12&2&122&122&1112\\\hline
${\mathcal E}_{1,j}-{\mathcal E}_{0,j}$&0   &0 &1 &0 &1&0  &1  &0   \\
${\mathcal E}_{2,j}-{\mathcal E}_{1,j}$&0   &0 &1 &0 &0&0  &0  &1   \\
${\mathcal E}_{3,j}-{\mathcal E}_{2,j}$&0   &0 &0 &1 &0&0  &0  &0   \\
${\mathcal E}_{4,j}-{\mathcal E}_{3,j}$&0   &0 &0 &0 &0&1  &0  &0   \\
${\mathcal E}_{5,j}-{\mathcal E}_{4,j}$&0   &0 &0 &0 &0&1  &0  &0   \\
${\mathcal E}_{6,j}-{\mathcal E}_{5,j}$&0   &0 &0 &0 &0&0  &1  &0   \\
${\mathcal E}_{7,j}-{\mathcal E}_{6,j}$&0   &0 &0 &0 &0&0  &0  &0   \\
\end{tabular}
\end{center}
Following Step (ii) and Step (iii), letters 1
contained in the above table are classified into 3 groups,
as indicated in the following table.
\begin{center}
\begin{tabular}{l|cccccccc}
                 &1111&11&22&12&2         &122&122       &1112\\\hline
${\mathcal E}_{1,j}-{\mathcal E}_{0,j}$&    &  &3 &  &$2^{\ast}$&   &1         &    \\
${\mathcal E}_{2,j}-{\mathcal E}_{1,j}$&    &  &3 &  &          &   &          &$1^{\ast}$\\
${\mathcal E}_{3,j}-{\mathcal E}_{2,j}$&    &  &  &3 &          &   &          &    \\
${\mathcal E}_{4,j}-{\mathcal E}_{3,j}$&    &  &  &  &          &3  &          &    \\
${\mathcal E}_{5,j}-{\mathcal E}_{4,j}$&    &  &  &  &          &3  &          &    \\
${\mathcal E}_{6,j}-{\mathcal E}_{5,j}$&    &  &  &  &          &   &$3^{\ast}$&    \\
${\mathcal E}_{7,j}-{\mathcal E}_{6,j}$&    &  &  &  &          &   &          &    \\
\end{tabular}
\end{center}
{}The cardinalities of the 3 groups are
2, 1 and 6, respectively.
{}From the positions marked with $\ast$, we find 
$(\mu_1,j_1)=(2,8)$,
$(\mu_2,j_2)=(1,5)$ and $(\mu_3,j_3)=(6,7)$.
Now we evaluate riggings $r_i$ according to the rule 
(\ref{s:eq:rigging}).
\begin{eqnarray*}
r_1&=&\sum_{i=1}^{8-1} \min(2,\lambda_i)
      +{\mathcal E}_{2,8}-2\sum_{i=1}^8{\mathcal E}_{2,i}\\
   &=&(2+2+2+2+1+2+2)+1-2(0+0+2+0+1+0+1+1)\\
   &=&4,\\
r_2&=&\sum_{i=1}^{5-1} \min(1,\lambda_i)
      +{\mathcal E}_{1,5}-2\sum_{i=1}^5{\mathcal E}_{1,i}\\
   &=&(1+1+1+1)+1-2(0+0+1+0+1)\\
   &=&1,\\
r_3&=&\sum_{i=1}^{7-1} \min(6,\lambda_i)
      +{\mathcal E}_{6,7}-2\sum_{i=1}^7{\mathcal E}_{6,i}\\
   &=&(4+2+2+2+1+3)+2-2(0+0+2+1+1+2+2)\\
   &=&0.
\end{eqnarray*}
Therefore we obtain
$(\mu_1,r_1)=(2,4)$, $(\mu_2,r_2)=(1,1)$ and $(\mu_3,r_3)=(6,0)$,
which coincide with the calculation in Example \ref{s:ex:kkr}.

The reader should
compare the above local energy distribution and
box adding procedure fully exhibited in Example \ref{s:ex:kkr}.
Then it will be observed that the complicated combinatorial procedure
in Definition \ref{s:def:kkr} is reduced 
to rather automatic applications of the combinatorial $R$
and energy functions.
\end{example}

\begin{example}\label{ex:big}
Theorem \ref{s:main} provides a panoramic view on the combinatorial
procedure of the KKR bijection 
from energy distribution.
To show a typical example, we pick 
the following long path
(length 30).
\begin{eqnarray*}
&&\fbox{22}\otimes\fbox{2}\otimes\fbox{2}\otimes\fbox{1}\otimes
\fbox{1122}\otimes\fbox{112}\otimes\fbox{1}\otimes\fbox{11}
\otimes\fbox{222}\otimes\fbox{12}\otimes\fbox{11}\otimes
\fbox{2}\\
&&\otimes\fbox{2}\otimes\fbox{2}\otimes\fbox{22}\otimes
\fbox{2}\otimes\fbox{1122}\otimes\fbox{22}\otimes\fbox{2}
\otimes\fbox{222}\otimes\fbox{1}\otimes\fbox{112}\otimes
\fbox{1}\otimes\fbox{12}\\
&&\otimes\fbox{1222}\otimes
\fbox{11122}\otimes\fbox{2}\otimes\fbox{22}\otimes\fbox{2}
\otimes\fbox{2}
\end{eqnarray*}
Then, the local energy distribution takes
the following form.
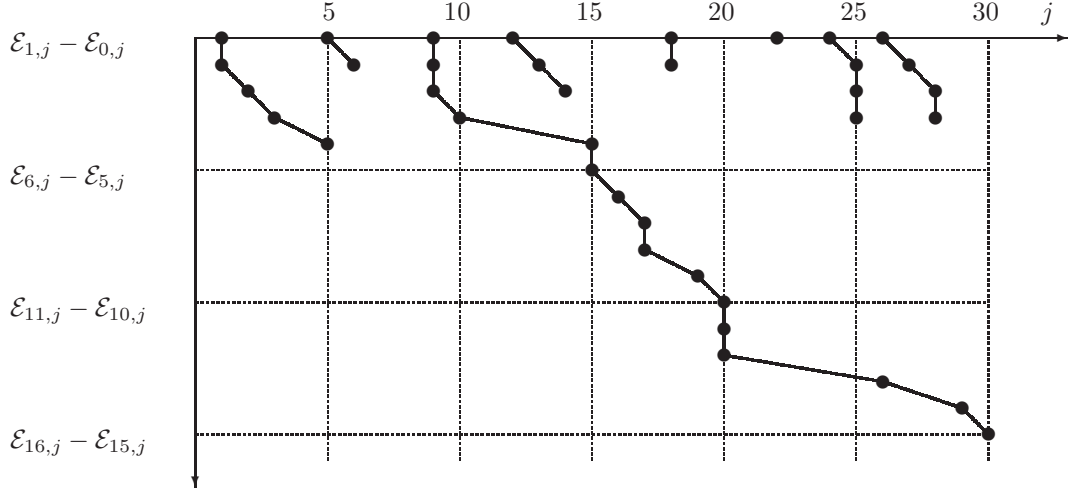
\begin{center}
\unitlength 10pt
\begin{picture}(39,18)(-6,-0.5)
\put(30,1){\circle*{0.5}}
\put(29,2){\circle*{0.5}}
\put(26,3){\circle*{0.5}}
\put(20,4){\circle*{0.5}}
\put(20,5){\circle*{0.5}}
\put(20,6){\circle*{0.5}}
\put(19,7){\circle*{0.5}}
\put(17,8){\circle*{0.5}}
\put(17,9){\circle*{0.5}}
\put(16,10){\circle*{0.5}}
\put(15,11){\circle*{0.5}}
\put(5,12){\circle*{0.5}}
\put(15,12){\circle*{0.5}}
\put(3,13){\circle*{0.5}}
\put(10,13){\circle*{0.5}}
\put(25,13){\circle*{0.5}}
\put(28,13){\circle*{0.5}}
\put(2,14){\circle*{0.5}}
\put(9,14){\circle*{0.5}}
\put(14,14){\circle*{0.5}}
\put(25,14){\circle*{0.5}}
\put(28,14){\circle*{0.5}}
\put(1,15){\circle*{0.5}}
\put(6,15){\circle*{0.5}}
\put(9,15){\circle*{0.5}}
\put(13,15){\circle*{0.5}}
\put(18,15){\circle*{0.5}}
\put(25,15){\circle*{0.5}}
\put(27,15){\circle*{0.5}}
\put(1,16){\circle*{0.5}}
\put(5,16){\circle*{0.5}}
\put(9,16){\circle*{0.5}}
\put(12,16){\circle*{0.5}}
\put(18,16){\circle*{0.5}}
\put(22,16){\circle*{0.5}}
\put(24,16){\circle*{0.5}}
\put(26,16){\circle*{0.5}}
\put(0,16){\vector(1,0){33}}
\put(0,16){\vector(0,-1){17}}
\multiput(5,0)(5,0){6}{
\multiput(0,0)(0,0.2){80}{\line(0,1){0.1}}
}
\multiput(0,1)(0,5){3}{
\multiput(0,0)(0.2,0){150}{\line(1,0){0.1}}
}
\put(4.8,16.7){5}
\put(9.4,16.7){10}
\put(14.4,16.7){15}
\put(19.4,16.7){20}
\put(24.4,16.7){25}
\put(29.4,16.7){30}
\put(32,16.7){$j$}
\put(-7,15.5){${\mathcal E}_{1,j}-{\mathcal E}_{0,j}$}
\put(-7,10.5){${\mathcal E}_{6,j}-{\mathcal E}_{5,j}$}
\put(-7,5.5){${\mathcal E}_{11,j}-{\mathcal E}_{10,j}$}
\put(-7,0.5){${\mathcal E}_{16,j}-{\mathcal E}_{15,j}$}
\thicklines
\qbezier(26,16)(26,16)(28,14)
\qbezier(28,14)(28,14)(28,13)
\qbezier(24,16)(24,16)(25,15)
\qbezier(25,15)(25,15)(25,13)
\qbezier(18,16)(18,16)(18,15)
\qbezier(12,16)(12,16)(14,14)
\qbezier(9,16)(9,16)(9,14)
\qbezier(9,14)(9,14)(10,13)
\qbezier(10,13)(10,13)(15,12)
\qbezier(15,12)(15,12)(15,11)
\qbezier(15,11)(15,11)(17,9)
\qbezier(17,9)(17,9)(17,8)
\qbezier(17,8)(17,8)(19,7)
\qbezier(19,7)(19,7)(20,6)
\qbezier(20,6)(20,6)(20,4)
\qbezier(20,4)(20,4)(26,3)
\qbezier(26,3)(26,3)(29,2)
\qbezier(29,2)(29,2)(30,1)
\qbezier(5,16)(5,16)(6,15)
\qbezier(1,16)(1,16)(1,15)
\qbezier(1,15)(1,15)(3,13)
\qbezier(3,13)(3,13)(5,12)
\end{picture}
\end{center}
In the above table, letters 1 in the local energy distribution
are represented by ``$\bullet$", and letters 0 are suppressed.
According to Step (ii) and Step (iii), $\bullet$ 
belonging to the same group are joined by thick lines.
We see there are 8 groups whose cardinalities are
5, 2, 16, 3, 2, 1, 4, 4 from left to right, respectively.

By using the formula (\ref{s:eq:rigging}), we get
the unrestricted rigged configuration as follows:
$(\mu_1,r_1)=(4,8)$, $(\mu_2,r_2)=(4,8)$,
$(\mu_3,r_3)=(1,10)$, $(\mu_4,r_4)=(2,8)$,
$(\mu_5,r_5)=(3,2)$, $(\mu_6,r_6)=(16,-15)$,
$(\mu_7,r_7)=(2,0)$, $(\mu_8,r_8)=(5,-5)$.
The vacancy numbers for each row is
$p_{16}=-15$, $p_5=7$,
$p_4=10$, $p_3=14$,
$p_2=16$ and $p_1=14$.
Note that since the path in this example is not
highest, the resulting unrestricted rigged configuration
has negative riggings and vacancy numbers.

\end{example}

\section{Generalized periodic box-ball system}\label{sec:pbbs}

Here we extend the inverse scattering formalism \cite{KTT} of the
simplest periodic box-ball system \cite{YT,YYT} to general higher spins.
The relevant time evolutions and associated energy 
will be denoted by $T_l$ and $E_l$ for distinction from 
${\mathcal T}_l$ and ${\mathcal E}_l$ for the non-periodic case.
Most of the proofs will be omitted as they are 
similar (but somewhat more involved) to \cite{KTT}.
Our new algorithm for the KKR map  
(Theorem \ref{s:main}), adapted to 
the periodic boundary condition, serves as a simple algorithm 
for the direct scattering transform.

\subsection{Time evolution}\label{subsec:te}
Fix the integer $L, s\in \Z_{\ge 1}$ throughout.
Set 
\begin{equation}\label{eq:pdef}
\Pth = B_s^{\otimes L}.
\end{equation}
We will also write $\Aff(\Pth) = \Aff(B_s)^{\otimes L}$.
An element of $\Pth$ is called a path. 
A path $b$ is {\em highest} if ${\tilde e}_1b = 0$.
The weight of a path $b=b_1\otimes \cdots \otimes b_L$ is 
given by 
$\wt(b) = \wt(b_1) + \cdots + \wt(b_L)$.
We write 
$\wt(b)>0$ ($\wt(b)<0)$ 
when it belongs to 
$\Z_{>0}\Lambda_1$ $(\Z_{<0}\Lambda_1)$.

Our generalized periodic box-ball system is a dynamical system 
on a subset of $\Pth$ equipped with the commuting family of 
time evolutions $T_1, T_2, \ldots$.
Let $b=b_1\otimes \cdots \otimes b_L \in \Pth$ be a path and 
$l \in \Z_{\ge 1}$.
For $v_l \in B_l$, suppose 
\begin{equation}\label{eq:bpe}
\zeta^0v_l \otimes (\zeta^0b_1 \otimes \cdots \otimes \zeta^0b_L) 
\simeq 
(\zeta^{-d_1}{\tilde b}_1 \otimes 
\cdots \otimes \zeta^{-d_L}{\tilde b}_L) \otimes 
\zeta^e v'_l
\end{equation}
holds under the isomorphism 
$\Aff(B_l)\otimes \Aff(\Pth) \simeq 
\Aff(\Pth) \otimes \Aff(B_l)$, 
where the right hand side is unambiguously determined from the 
left hand side. ($e=d_1+\cdots +d_L$.)
We say that $b$ is {\em $T_l$-evolvable} if 
the following (i) existence and (ii) uniqueness are satisfied:

\begin{enumerate}
\item  there exists $v_l \in B_l$ such that $v'_l=v_l$.

\item if there are more than one such $v_l$, 
${\tilde b}_1 \otimes \cdots \otimes {\tilde b}_L$ is 
independent of their choice.
\end{enumerate}
If $b$ is $T_l$-evolvable,  we define 
$T_l(b) = {\tilde b}_1 \ot \cdots \ot {\tilde b}_L (\in \Pth)$.
Otherwise we set $T_l(b) = 0$.
In this sense, we will also write $T_l(b) \neq 0$ to mean 
that $b$ is $T_l$-evolvable.
\begin{lemma}\label{le:eu}
If $b=b_1\otimes \cdots \otimes b_L$ is $T_l$-evolvable, not only 
${\tilde b}_1 \otimes \cdots \otimes {\tilde b}_L$ but also 
$d_1, \ldots, d_L$ and $e$ in (\ref{eq:bpe}) are independent of the 
possibly nonunique choices of $v_l=v'_l$.
\end{lemma}
Thanks to this lemma we are entitled to define
$E_l(b) =  e ( \in \Z_{\ge 0})$ by (\ref{eq:bpe}) for a
$T_l$-evolvable path $b$.
Actually $v_l$ can be nonunique only if $l > s$ and 
$\wt(p)=0$.
The operations $T_1, T_2, \ldots$  form a family of time evolution operators 
associated with the energy $E_1, E_2, \ldots$.
These definitions can be summarized in
\begin{equation}\label{eq:tle}
\zeta^0 v_l \otimes b \simeq T_l(b) \otimes \zeta^{E_l(b)}v_l
\end{equation} 
up to the spectral parameter for $T_l$-evolvable $b$.
Pictorially, (\ref{eq:bpe}) looks as
\begin{equation}\label{eq:vb}
\begin{picture}(200,35)(10,-1)

\put(-30,13){$v_l=v(0)$}

\multiput(0,0)(40,0){1}{
\put(10,15){\line(1,0){20}}\put(20,9){\line(0,1){10}}
\put(17,22){$b_1$}\put(17,-4){${\tilde b}_1$}
\put(34,13){$v(1)$}}

\multiput(45,0)(40,0){1}{
\put(10,15){\line(1,0){20}}\put(20,9){\line(0,1){10}}
\put(17,22){$b_2$}\put(17,-4){${\tilde b}_2$}
\put(34,13){$v(2)$}}

\put(104,13){$\cdots$}

\multiput(110,0)(40,0){1}{
\put(10,15){\line(1,0){20}}\put(20,9){\line(0,1){10}}
\put(18,22){$b_{L-1}$}\put(17,-4){${\tilde b}_{L-1}$}
\put(34,13){$v(L-1)$}}

\multiput(175,0)(40,0){1}{
\put(10,15){\line(1,0){20}}\put(20,9){\line(0,1){10}}
\put(18,22){$b_{L}$}\put(17,-4){${\tilde b}_{L}$}
\put(34,13){$v(L)=v'_l$}}

\end{picture}
\end{equation} 
Clearly the time evolutions are 
weight preserving, i.e., $\wt(T_l(b)) = \wt(b)$ when 
$T_l(b) \neq 0$.

Since the combinatorial $R$ is trivial 
on $B_s \ot B_s$ (see (\ref{eq:tri})),
we have the unique choice $v_s=v'_s=b_L$ in (\ref{eq:bpe}),
saying that a path is always $T_s$-evolvable and 
$T_s$ acts as a cyclic shift: 
\begin{equation}\label{eq:ts}
T_s(b_1 \otimes b_2 \otimes \cdots \otimes b_L)
= b_L \otimes b_1 \otimes \cdots \otimes b_{L-1}.
\end{equation}
If $s=1$, all the paths are $T_l$-evolvable 
for any $l \ge 1$ \cite{KTT}. 
However this is no longer the case for $s>1$. 
A similar situation is known also in the higher rank extensions
\cite{KT}.
Here we treat such a subtlety 
characteristic in the periodic setting.

We simply say that $b \in \Pth$ is {\em evolvable}
if it is $T_l$-evolvable for all $l \in \Z_{\ge 1}$.
We warn that ``$b$ is $T_l$-evolvable" is different from 
``$T_l(b)$ is evolvable". 
The former means $T_l(b) \neq 0$ whereas the latter 
does $T_kT_l(b) \neq 0$ for all $k \ge 1$.
Here is a characterization of evolvable paths.
\begin{proposition}\label{pr:evo}
A path $b =b_1 \otimes \cdots \otimes b_L$ is 
evolvable if and only if 
$b_i=\fbox{1...1}$ or $b_i=\fbox{2...2}$ for some $i$.
\end{proposition}
The proof of the proposition also tells the way  
to construct $v_l$ that makes (\ref{eq:tle}) hold for 
a given path $b$.
For $l \ge s$, determine $v_l \in B_l$ by
(see (\ref{eq:ul}) for $u_l$)
\begin{equation}\label{eq:vt}
\begin{split}
u_l \otimes b \simeq b' &\otimes v_l 
\quad \hbox{ if } \wt(b)\ge 0,\\
\omega(u_l) \otimes b \simeq b' &\otimes v_l 
\quad \hbox{ if } \wt(b) < 0,
\end{split}
\end{equation}
where $b' \in \Pth$ is another path.
So obtained $v_l$ is shown to satisfy 
$v_l \ot b \simeq T_l(b)\ot v_l$ under
$B_l \ot \Pth \simeq \Pth \ot B_l$.
One may either use the latter relation 
in (\ref{eq:vt}) to define $v_l$ when $\wt(b)=0$.
For $l<s$, one has the unique $v_l$ from 
$v \ot b \simeq b' \ot v_l$ for arbitrary $v \in B_l$ if 
$b$ is evolvable. 
Then $v_l \ot b \simeq T_l(b)\ot v_l$ is again valid.

\begin{theorem}\label{th:te}
Suppose $b, T_l(b)$ and $T_k(b)$ are evolvable.
Then the commutativity $T_lT_k(b) = T_kT_l(b)$
and the energy 
conservation $E_l(T_k(b))=E_l(b), E_k(T_l(b))=E_k(b)$ hold.
\end{theorem}
\begin{proof}
Take $v_k$ for $b$ and $v_l$ for $T_k(b)$ as in (\ref{eq:vt}).
Set $R(\zeta^0v_l \otimes \zeta^0 v_k) = 
\zeta^{-\delta}\overline{v}_k \otimes 
\zeta^{\delta}\overline{v}_l$ and 
regard $b$ as an element of $\Aff(\Pth)$.
By using the combinatorial $R$, one can reorder 
$\zeta^0v_l \otimes \zeta^0v_k \otimes b$ 
in two ways along the isomorphism 
$\Aff(B_l) \otimes \Aff(B_k) \otimes 
\Aff(\Pth) \simeq 
\Aff(\Pth) \otimes \Aff(B_k) \otimes \Aff(B_l)$
as follows:

\setlength{\unitlength}{1mm}
\begin{picture}(80,55)(-20,0)

\put(0,10){\line(1,0){8}}
\put(0,20){\line(1,0){8}}

\put(11,9){$\cdots$} 
\put(12,19){$\cdots$}

\put(18,10){\line(1,0){8}}
\put(18,20){\line(1,0){8}}

\put(4,7){\line(0,1){16}}
\put(22,7){\line(0,1){16}}

\put(-4,19){$\scriptstyle v_{k}$}
\put(-4,9){$\scriptstyle v_{l}$}

\put(12,25){$\scriptstyle b$}
\put(11,14){$\scriptstyle T_{k}(b)$}
\put(10,4){$\scriptstyle T_lT_{k}(b)$}

\put(28,19){$\scriptstyle \zeta^{E_{k}(b)}{v}_{k}$}
\put(27,9){$\scriptstyle \zeta^{E_{l}(T_{k}(b))}{v}_{l}$}

\put(44,10){\line(1,1){10}}
\put(44,20){\line(1,-1){10}}

\put(55,9){$\scriptstyle \zeta^{E_{k}(b)-\delta}\overline{v}_{k}$}
\put(55,19){$\scriptstyle \zeta^{E_{l}(T_{k}(b))+\delta}\overline{v}_{l}$}

\put(-10,14){$=$}

\put(0.2,35){\line(1,1){10}}
\put(0.2,45){\line(1,-1){10}}

\put(-4,35){$\scriptstyle v_{l}$}
\put(-4,45){$\scriptstyle v_{k}$}

\put(11,35){$\scriptstyle \zeta^{-\delta}\overline{v}_{k}$}
\put(12,45){$\scriptstyle \zeta^{\delta}\overline{v}_{l}$}

\put(20,35){\line(1,0){8}}
\put(20,45){\line(1,0){8}}

\put(31,34){$\cdots$} 
\put(32,44){$\cdots$}

\put(38,35){\line(1,0){8}}
\put(38,45){\line(1,0){8}}

\put(24,32){\line(0,1){16}}
\put(42,32){\line(0,1){16}}

\put(-4,19){$\scriptstyle v_{k}$}
\put(-4,9){$\scriptstyle v_{l}$}

\put(32,50){$\scriptstyle b$}
\put(31,39){$\scriptstyle T_{l}(b)$}
\put(30,29){$\scriptstyle T_{k}T_{l}(b)$}

\put(47,45){$\scriptstyle \zeta^{E_{l}(b)+\delta}\overline{v}_{l}$}
\put(47,35){$\scriptstyle \zeta^{E_{k}(T_{l}(b))-\delta}\overline{v}_{k}$}

\put(75,9){,}

\end{picture}

\noindent
where the equality is due to the 
Yang-Baxter equation.
The outputs have been identified  
with $T_kT_l(b), \zeta^{E_{k}(T_{l}(b))-\delta}\overline{v}_{k}$, etc.
In particular the uniqueness (ii) stated under (\ref{eq:bpe}) 
guarantees that $\overline{v}_k\otimes T_l(b) \simeq 
T_kT_l(b) \otimes \overline{v}_k$ and 
$\overline{v}_l\otimes b\simeq 
T_l(b) \otimes \overline{v}_l$ up to the spectral parameter.
The sought relations 
$T_lT_k(b) = T_kT_l(b)$ and 
$E_l(T_k(b))=E_l(b), E_k(T_l(b))=E_k(b)$
are obtained by comparing the two sides.
\end{proof}

Let $s_0,\, s_1$ be the Weyl group operators (\ref{eq:wdef})
and $\omega$ be the involution (\ref{eq:odef}) 
acting on the crystal $\Pth$.
Then ${\widehat W}(A^{(1)}_1)
= \langle \omega, s_0, s_1\rangle$ forms 
the extended affine Weyl group of type $A^{(1)}_1$.
The time evolutions $T_l$ and the energy $E_l$ 
enjoy the symmetry under 
${\widehat W}(A^{(1)}_1)$.
\begin{proposition}\label{pr:w}
Let $b$ be an evolvable path.
Then for any $w \in {\widehat W}(A^{(1)}_1)$, 
$w(b)$ is also evolvable and 
the commutativity 
$wT_l(b) = T_l(w(b))$ and 
the invariance $E_l(w(b)) = E_l(b)$ are valid.
\end{proposition}
In particular, the relation
\begin{equation}\label{eq:tw}
T_l = \omega \circ T_l \circ \omega
\end{equation}
exchanging the letters $1 \leftrightarrow 2$ 
is useful. 

Any path is $T_l$-evolvable for $l\ge s$.
In fact, for $l$ sufficiently large 
the time evolution $T_l$ and the energy $E_l$ 
admit a simple description as follows.
\begin{proposition}\label{pr:tws} 
For any path $b \in \Pth$, there exists $k\ge s$ 
such that $T_l(b)$ and $E_l(b)$ are independent of $l$ for $l \ge k$.
Denoting them by $T_\infty(b)$ and $E_\infty(b)$, one has
\begin{equation*}
\begin{split}
T_\infty(b)=\omega(s_0(b)),  \;\;
\wt(b) = p_\infty\Lambda_1 \;\;\;& \;\hbox{ if } \wt(b)\ge 0,\\
T_\infty(b)=\omega(s_1(b)),  \;\;
\wt(b) = -p_\infty\Lambda_1 & \;\hbox{ if } \wt(b) \le 0,
\end{split}
\end{equation*} 
where $p_\infty = Ls - 2E_\infty(b)$ according to (\ref{eq:va}).
In particular, $T_\infty(b)=\omega(b)$ if $\wt(b)=0$. 
\end{proposition}
\begin{example}\label{ex:ts}
For $b=\fbox{112}\otimes \fbox{111} \otimes \fbox{222}
\otimes \fbox{122} \otimes \fbox{112}$
having a positive weight, we have
\begin{align*}
T_1(b)&= \fbox{122}\otimes \fbox{111}\otimes \fbox{122}\otimes \fbox{222}\otimes \fbox{111},\\
T_2(b)&=\fbox{112}\otimes \fbox{112}\otimes \fbox{112}\otimes \fbox{222}\otimes \fbox{112},\\
T_3(b)&=\fbox{112}\otimes \fbox{112}\otimes \fbox{111}\otimes \fbox{222}\otimes \fbox{122},\\
T_4(b)&=\fbox{122}\otimes \fbox{112}\otimes \fbox{111}\otimes \fbox{122}\otimes \fbox{122},\\
T_l(b)&=\fbox{122}\otimes \fbox{122}\otimes \fbox{111}\otimes \fbox{112}\otimes \fbox{122} 
\;\; (l \ge 5). 
\end{align*}
So $T_\infty(b) = T_5(b)$. 
On the other hand, $0$-signature and reduced $0$-signature of $b$ read
\begin{equation*}
\underset{--}{\underset{--+}{\fbox{112}}}\otimes 
\underset{--}{\underset{---}{\fbox{111}}} \otimes 
\underset{++}{\underset{+++}{\fbox{222}}}
\otimes \underset{-++}{\fbox{122}} \otimes 
\underset{+}{\underset{--+}{\fbox{112}}}
\end{equation*}
Thus $s_0(b)=\fbox{112}\otimes \et{0}\fbox{111} \otimes \fbox{222}
\otimes \fbox{122} \otimes \fbox{112}$, which coincides with 
$\omega(T_\infty(b))$.
\end{example}

For an evolvable path $b \in \Pth$, we have the time evolution
$T_l(b) \in \Pth$ and the 
associated energy $E_l(b) \in \Z_{\ge 0}$ for all $l \ge 1$.
This leads us to introduce the ``iso-level" set 
\begin{equation}\label{eq:pm}
{\widehat \Pth}(m) = \{ b \in \Pth \mid b: \hbox{evolvable}, \; 
E_l(b) = \sum_{k\ge 1} \min(l,k)m_k\}
\end{equation}
labeled with the sequence $m=\{m_k\mid k \ge 1\}$.
We shall always take it for granted that 
$\{m_k\}$ and $\{E_l\}$ are in one-to-one correspondence via
\begin{equation}\label{eq:me}
E_l = \sum_{k\ge 1}\min(l,k)m_k,\quad 
m_k = -E_{k-1}+2E_k-E_{k+1}\; (E_0=0).
\end{equation}
We also use the vacancy number 
\begin{equation}\label{eq:va}
p_j = L\min(s,j) - 2E_j.
\end{equation}
The following result is due to T. Takagi.
\begin{proposition}[\cite{Ta}]\label{le:ta}
For any path $b \in {\widehat \Pth}(m)$ with $\wt(b)\ge 0$,
its time evolution $\bigl(\prod_l T_{l}^{d_l})(b)$ 
becomes highest under appropriate choices of $\{d_l\}$
\footnote{Actually $T_s$ and $T_{s-1}$ have been shown to suffice. }.
\end{proposition}
Such $\{d_l\}$ is not unique.
Cyclic shift $T_s^{d_s}$ is not enough to achieve this in general.
{}From Proposition \ref{le:ta} one can show
\begin{proposition}\label{pr:nonemp}
${\widehat \Pth}(m) \neq \emptyset$ if and only if 
$\forall p_j \ge 0$.
\end{proposition}
Henceforth we assume $\forall p_j \ge 0$.
If $b$ belongs to ${\widehat \Pth}(m)$ and 
$T_l(b)$ is evolvable, then $T_l(b) \in {\widehat \Pth}(m)$
must hold because of $E_k(T_l(b)) = E_k(b)$ by Theorem \ref{th:te}.
However, the point here is that 
even if a path $b$ is evolvable, 
it is {\em not} guaranteed in general that 
its time evolution $T_l(b)$ is again evolvable.
\begin{example}\label{ex:enc}
$b_1=\fbox{11}\ot \fbox{22}$ and 
$b_2=\fbox{22}\ot \fbox{11}$ are evolvable, but
$T_1(b_1) = T_1(b_2) = \fbox{12} \ot \fbox{12}$ is not.
See Proposition \ref{pr:evo}.
The situation is depicted as 
\begin{equation*}
\begin{picture}(150,50)(0,0)
\put(10,37){$\fbox{11}\otimes \fbox{22}$}
\put(60,37){\vector(2,-1){15}}\put(66,37){$T_1$}

\put(10,10){$\fbox{22}\otimes \fbox{11}$}
\put(60,17){\vector(2,1){15}}\put(67,10){$T_1$}

\put(80,25){$\fbox{12}\otimes \fbox{12}$}
\put(130,30){\vector(1,0){20}}\put(135,34){$T_1$}

\put(155,28){$0$}

\end{picture}
\end{equation*}
\end{example}
Thus the set $T_l({\widehat \Pth}(m))$ can contain non-evolvable paths
in general.
On the other hand, all the evolvable paths in 
$T_l({\widehat \Pth}(m))$ must share the same 
energy spectrum $\{E_l\}$ as ${\widehat \Pth}(m)$
by virtue of Theorem \ref{th:te}.
Therefore, what holds in general is  
\begin{equation*}
T_l({\widehat \Pth}(m)) = (\hbox{subset of } 
{\widehat \Pth}(m)) \sqcup 
\{T_l(b): \hbox{non-evolvable} \mid b \in {\widehat \Pth}(m)\}.
\end{equation*}
A natural question is to find a pleasant situation where 
$T_l$ acts on ${\widehat \Pth}(m)$ as a bijection.
This is answered in 
\begin{proposition}\label{pr:cyc}
$T_l({\widehat \Pth}(m)) = {\widehat \Pth}(m)$
holds for all $l$
if and only if $(E_1, E_2) \neq (L/2, L)$.
\end{proposition} 
So this is always satisfied if $L$ is odd.
For an evolvable path $b$ with even length $L$, 
the condition $(E_1, E_2) = (L/2, L)$ is 
equivalent to 
\begin{equation*}
b = 11(c_1)\otimes (c_2)22 \otimes \cdots \otimes (c_L)22\;\;
\hbox{ or } \;\;
b = (c_1)22\otimes 11(c_2) \otimes \cdots \otimes 11(c_L)
\end{equation*}
for some $c_i \in B_{s-2}$, where $11$ and $22$ alternate.
Here for example, $11(c)=\fbox{11122}$ and 
$(c)22=\fbox{12222}$ for $c=\fbox{122} \in B_3$.
(Thus such $b$ can exist only for $s\ge 2$.)
The two paths in Example \ref{ex:enc} correspond to the case
$(E_1, E_2) = (1,2)$ with $L=2$.
For ${\widehat \Pth}(m)$ such that $(E_1, E_2) \neq (L/2, L)$,
the inverse time evolution is given by 
\begin{equation}\label{eq:inv}
T^{-1}_l = \varrho \circ T_l \circ \varrho,
\end{equation}
where $\varrho$ is defined by 
\begin{equation*}
\varrho(b_1 \otimes b_2 \otimes  \cdots \otimes b_L) 
= b_L \otimes \cdots \otimes b_2 \otimes b_1.
\end{equation*}

To summarize Theorem \ref{th:te}, Proposition \ref{pr:w} and 
Proposition \ref{pr:cyc}, 
each set ${\widehat \Pth}(m)$ of evolvable paths is 
characterized by the 
conserved quantity $E_l = \sum_{k \ge 1}\min_k(l,k)m_k$ called 
energy, and enjoys the invariance under 
\begin{enumerate}
\item the extended affine Weyl group ${\widehat W}(A^{(1)}_1)$,
  
\item the commuting family of 
invertible time evolutions $\{T_l \mid l \ge 1\}$.
\end{enumerate}
${\widehat \Pth}(m)$ is non-empty if $\forall p_j \ge 0$.
The invariance (ii) is valid if 
$(E_1, E_2) \neq (L/2,L)$ is further satisfied.

\subsection{Action-angle variable}\label{ssec:aav}

{}From now on, we assume that $m=\{m_j\}$ satisfies
\begin{equation}\label{eq:pos}
\forall p_j \ge 1.
\end{equation}
See (\ref{eq:me}) and (\ref{eq:va}).
This fulfills the conditions in Propositions 
\ref{pr:nonemp} and \ref{pr:cyc} since $p_1 = L-2E_1$.
The set ${\widehat \Pth}(m)$ 
is decomposed into a disjoint union of fixed weight subsets:
\begin{equation}\label{eq:du}
\begin{split}
{\widehat \Pth}(m) &= \Pth(m) \sqcup \omega(\Pth(m)),\\
\Pth(m) &= 
\{b \in {\widehat \Pth}(m) \mid \wt(b) = p_\infty\Lambda_1 \}.
\end{split}
\end{equation}
In view of (\ref{eq:tw}), 
dynamics on ${\widehat \Pth}(m)$ is reduced to 
the commuting family of invertible time evolutions 
$\{T_l\}$ on the fixed (positive) weight subset $\Pth(m)$:
\begin{equation}\label{eq:tp}
T_l \, : \; \Pth(m) \rightarrow \Pth(m).
\end{equation}
We present the inverse scattering transform that linearizes
the dynamics (\ref{eq:tp}) and an explicit solution of the 
initial value problem.
For a general background on the inverse 
scattering method, see \cite{AS,GGKM}.
In our approach the direct scattering transform is 
formulated either by a modified KKR bijection 
as in the $s=1$ case \cite{KTT} or by an appropriate 
extension of the procedure in Theorem \ref{s:main} 
to a periodic setting.

First we introduce the action-angle variables.
For the paths belonging to $\Pth(m)$, 
the action variable is just the conserved quantity
$m=\{m_j\}$ or equivalently $\{E_l\}$ (\ref{eq:me}).
It may also be presented as the Young diagram 
$\mu$ having $m_j$ rows with length $j$.
Let $H = \{j_1 < \cdots < j_g\}$ be the 
set of distinct row lengths of $\mu$, namely,
$j \in H \leftrightarrow m_j>0$.
The set ${\mathcal J}(m)$ of angle variables is defined by
\begin{equation}\label{eq:av}
\begin{split}
{\mathcal J}(m) &= \Bigl(
(\Z^{m_{j_1}}\times 
\cdots \times \Z^{m_{j_g}})/\Gamma\;
-\Delta\Big)_{\rm sym},\\
\Gamma &= A(\Z^{m_{j_1}}\times 
\cdots \times \Z^{m_{j_g}}).
\end{split}
\end{equation}
Here $A=(A_{j\alpha, k\beta})$ is the matrix of size 
$m_{j_1} + \cdots + m_{j_g}$ having a block structure:
\begin{equation}\label{eq:a}
A_{j\alpha, k\beta}=\delta_{j,k}\delta_{\alpha,\beta}
(p_j+m_j) + 2\min(j,k)-\delta_{j,k},
\end{equation}
where $j,k \in H$ and $1\le \alpha \le m_j, 1 \le \beta \le m_k$.
$A$ is symmetric and positive definite under the assumption 
(\ref{eq:pos}) \cite{KN}.
$\Delta$ is the subset of 
$(\Z^{m_{j_1}}\times 
\cdots \times \Z^{m_{j_g}})/\Gamma$
having coincident components within a block:
\begin{equation}
\Delta = \{(I_{j,\alpha})_{j \in H, 1 \le \alpha \le m_j}
\mid I_{j,\alpha}=I_{j,\beta} \hbox{ for some }
j \in H, 1 \le \alpha \neq \beta \le m_j\}.
\end{equation}
In (\ref{eq:av}), 
$-\Delta$ signifies the complement of $\Delta$.
The subscript ${\rm sym}$ 
means the identification under the exchange of 
components within blocks via the symmetric group
${\mathfrak S}_{m_{j_1}} \times \cdots \times 
{\mathfrak S}_{m_{j_g}}$.
We introduce the time evolution of the angle variables by
\begin{equation*}
\begin{split}
T_l\, : \; {\mathcal J}(m) &\longrightarrow {\mathcal J}(m),\\
(I_{j,\alpha}) & \longmapsto 
(I_{j,\alpha} + \min(l,j)),
\end{split}
\end{equation*}
which makes sense because it 
obviously preserves the set 
$(\Z^{m_{j_1}}\times 
\cdots \times \Z^{m_{j_g}})/\Gamma\;-\Delta$.
We shall simply write this as
\begin{equation}\label{eq:tla}
T_l({\bf I}) = {\bf I} + {\bf h}_l.
\end{equation}
Namely,  
${\bf h}_l = (\min(j,l))_{j \in H, 1 \le \alpha \le m_j}
\in \Z^{m_{j_1}+\cdots +m_{j_g}}$ 
is the velocity of the angle variable ${\bf I} = (I_{j,\alpha})$
under the time evolution $T_l$.

\subsection{Direct scattering}\label{ssec:dst}
We introduce the direct scattering map
$\Phi: \Pth(m) \rightarrow {\mathcal J}(m)$.
A quick formulation is due to a modified KKR bijection as 
done in \cite{KTT} for $s=1$. 
Note that Proposition \ref{le:ta} tells that 
${\mathcal P}(m)$ is the $\{T_l\}$-orbit of highest paths 
having the KKR configuration $m=\{m_j\}$.
Thus we express a given $b \in \Pth(m)$ as 
$b = \bigl(\prod_l T_{l}^{d_l})(b_+)$ in terms of a highest path 
$b_+ \in \Pth(m)$.
Let $(\mu, r)$ be the rigged configuration for $b_+$,
where the appearance of the $\mu$ corresponding to 
$m=\{m_j\}$ is due to Theorem \ref{th:te}.
Let the rigging attached to the length $j (\in H)$ rows of $\mu$
be $0 \le r_{j,1} \le \cdots \le r_{j,m_j} \le p_j$.
Consider the element 
\begin{equation}
{\bf I} + \sum_l d_l{\bf h}_l\;\; \hbox{\rm mod}\,\;\Gamma \;\;
\in {\mathcal J}(m),\;\;
\text{where}\;\;
{\bf I} = (r_{j,\alpha}+\alpha-1)_{j \in H, 1\le \alpha \le m_j}.
\label{eq:Phi}
\end{equation}
$r_{j,\alpha}+\alpha-1$ is strictly increasing with 
$\alpha$, therefore 
${\bf I} + \sum_l d_l{\bf h}_l\, \hbox{ mod } \Gamma$ 
belongs to 
$(\Z^{m_{j_1}}\times 
\cdots \times \Z^{m_{j_g}})/\Gamma-\Delta$.
Given $b \in \Pth(m)$, the choice of $\{d_l\}$ and 
the highest path $b_+$ such that  
$b=\bigl(\prod_l T_{l}^{d_l})(b_+)$ is not unique in general.
This non-uniqueness is to be cancelled by 
$\hbox{\rm mod}\, \Gamma$.
In fact we have 
\begin{theorem}\label{th:cd}
The rule $\Phi(b) = {\bf I} + \sum_l d_l{\bf h}_l\; \hbox{\rm mod}\, 
\Gamma$
specified by (\ref{eq:Phi}) is a bijection 
$\Phi: \Pth(m) \rightarrow {\mathcal J}(m)$, and 
the following commutative diagram is valid:
\begin{equation}\label{eq:cd}
\begin{CD}
{\mathcal P}(m) @>{\Phi}>> \J(m) \\
@V{T_l}VV @VV{T_l}V\\
{\mathcal P}(m) @>{\Phi}>> \J(m) 
\end{CD}
\end{equation}
\end{theorem}

An alternative way to define the 
direct scattering map $\Phi$ is obtained by a periodic extension 
of the procedure (i), (ii)', (iii), (iv) in Theorem \ref{s:main}.
This option is valid under a certain condition 
which we shall explain after (\ref{eq:av1}).
It is more direct than the above one in that the relation 
$b=\bigl(\prod_l T_{l}^{d_l})(b_+)$ need not be found.
Here we illustrate it along the example:
\begin{equation}\label{eq:bex}
b = \fbox{122} \otimes \fbox{122} \otimes \fbox{112} \otimes 
\fbox{112} \otimes \fbox{111} \otimes \fbox{122} \otimes 
\fbox{111} \otimes \fbox{111} \otimes \fbox{112}
\in B_3^{\otimes 9}\, .
\end{equation} 
This path has appeared as the $t=0$ case in the introduction.
Let $v_l \in B_l$ be the element satisfying 
$v_l\otimes b \simeq T_l(b)\otimes v_l$.
It is unique under the condition (\ref{eq:pos}) and 
can be found by (\ref{eq:vt}).
The energy is given by 
$E_1 = 4, E_2 = 7, E_3 = 8, E_l = 9 \,(l\ge 4)$.
So the action variable is 
\begin{equation*}
\mu = \yng(4,2,2,1)
\end{equation*}
Local energy $E_{l,k}=H(v(k-1)\otimes b_k)$ is determined 
by using $v(k)$ in (\ref{eq:vb}).
The distribution of $\delta E_{l,k} = E_{l,k}- E_{l-1,k}$ looks as
\begin{center}
\begin{tabular}{l|ccccccccc}
& 122 & 122 & 112 & 112 & 111 & 122 
& 111 & 111 & 112 \\\hline
$\delta E_{1,k}$&0   &1 &0 &1 &0&1  &0  &0   &1   \\
$\delta E_{2,k}$&1   &0 &1 &0 &0&1  &0  &0   &0\\
$\delta E_{3,k}$&1   &0 &0 &0 &0&0  &0  &0   &0\\
$\delta E_{4,k}$&0   &1 &0 &0 &0&0  &0  &0   &0\\
\end{tabular}
\end{center}
We group $1$'s by a periodic analogue of the 
procedure (i), (ii)', (iii), (iv) in Theorem \ref{s:main}.
Pick a lowest $1$, say $\delta E_{l,k}=1$ at the $l$th row.
If there are more than one such $k$, 
any choice is possible.
Let the rightmost $1$ in 
\begin{equation}
\delta E_{l-1,k+1},\ldots,
\delta E_{l-1,L-1}, \delta E_{l-1,L}, \delta E_{l-1,1},\ldots, 
\delta E_{l-1,k-1}, \delta E_{l-1,k}
\end{equation}
be $\delta E_{l-1,k'}=1$. Namely, $k'$ is the position of the 
rightmost $1$ satisfying $k' \le k$ cyclically.
Then connect $\delta E_{l,k}$ to $\delta E_{l-1,k'}$.
Repeat this until the successive connection 
reaches some $\delta E_{1,k''}$ on the first row.
This completes one group. Erase all the $1$'s in it and 
repeat the same procedure starting from a lowest $1$
in the rest 
to form other groups until all the initial $1$'s are exhausted.

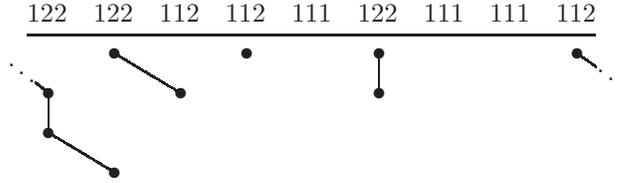
\begin{center}
\unitlength 10pt
\begin{picture}(24,8)(-0.2,4)
\put(0,10.2){$
122 \;\;\;\, 122 \;\;\;\, 112 \;\;\;\, 112 \;\;\;\, 111 \;\;\;\, 
122 \;\;\;\, 111 \;\;\;\, 111 \;\;\;\, 112$}
\put(0,9.7){\line(1,0){21.5}}

\put(3.3,9){\circle*{0.4}}
\put(8.3,9){\circle*{0.4}}
\put(13.3,9){\circle*{0.4}}
\put(20.8,9){\circle*{0.4}}

\put(0.8,7.5){\circle*{0.4}}
\put(5.8,7.5){\circle*{0.4}}
\put(13.3,7.5){\circle*{0.4}}

\put(0.8,6){\circle*{0.4}}

\put(3.3,4.5){\circle*{0.4}}

\qbezier(13.3,7.5)(13.3,7.5)(13.3,9)

\qbezier(5.8,7.5)(3.3,9)(3.3,9)

\qbezier(20.8,9)(20.8,9)(21.5,8.5)
\put(21.5,7.6){$\ddots$}

\qbezier(0.8,7.5)(0.8,7.5)(0.3,7.9)
\put(-0.8,7.8){$\ddots$}
\qbezier(0.8,6)(0.8,6)(0.8,7.5)
\qbezier(3.3,4.5)(3.3,4.5)(0.8,6)

\end{picture}
\end{center}

\noindent
A group consisting of $l$ dots will be called a soliton of length $l$.
Make a cyclic shift $T_s^{-d}$ so that all the solitons
stay within the left and the right boundary.
Namely, no soliton sits across the boundary.
In the above, we take for example $d=3$.

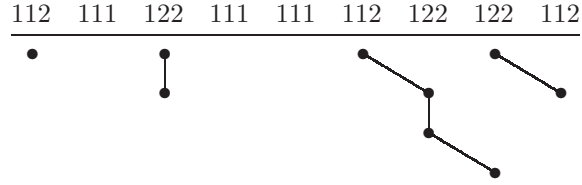
\begin{center}
\unitlength 10pt
\begin{picture}(24,8)(-0.2,4)
\put(0,10.2){$
112 \;\;\;\, 111 \;\;\;\, 
122 \;\;\;\, 111 \;\;\;\, 111 \;\;\;\, 112  \;\;\;\, 
122 \;\;\;\, 122 \;\;\;\, 112$}
\put(0,9.7){\line(1,0){21.5}}

\put(0.8,9){\circle*{0.4}}
\put(5.8,9){\circle*{0.4}}
\put(13.3,9){\circle*{0.4}}
\put(18.3,9){\circle*{0.4}}

\put(5.8,7.5){\circle*{0.4}}
\put(15.8,7.5){\circle*{0.4}}
\put(20.8,7.5){\circle*{0.4}}

\put(15.8,6){\circle*{0.4}}

\put(18.3,4.5){\circle*{0.4}}

\qbezier(5.8,9)(5.8,9)(5.8,7.5)
\qbezier(13.3,9)(13.3,9)(15.8,7.5)
\qbezier(15.8,7.5)(15.8,6)(15.8,6)
\qbezier(18.3,9)(18.3,9)(20.8,7.5)
\qbezier(15.8,6)(15.8,6)(18.3,4.5)

\end{picture}
\end{center}

\noindent
Computing the rigging of each soliton according to 
(\ref{s:eq:rigging}), we find 

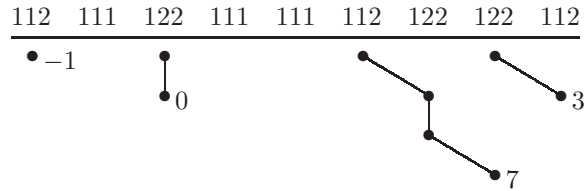
\begin{center}
\unitlength 10pt
\begin{picture}(24,8)(-0.2,4)
\put(0,10.2){$
112 \;\;\;\, 111 \;\;\;\, 
122 \;\;\;\, 111 \;\;\;\, 111 \;\;\;\, 112  \;\;\;\, 
122 \;\;\;\, 122 \;\;\;\, 112$}
\put(0,9.7){\line(1,0){21.5}}

\put(0.8,9){\circle*{0.4}}\put(1.2,8.5){$-1$}
\put(5.8,9){\circle*{0.4}}
\put(13.3,9){\circle*{0.4}}
\put(18.3,9){\circle*{0.4}}

\put(5.8,7.5){\circle*{0.4}}\put(6.2,7){$0$}
\put(15.8,7.5){\circle*{0.4}}
\put(20.8,7.5){\circle*{0.4}}\put(21.2,7){$3$}

\put(15.8,6){\circle*{0.4}}

\put(18.3,4.5){\circle*{0.4}}\put(18.7,4){$7$}

\qbezier(5.8,9)(5.8,9)(5.8,7.5)
\qbezier(13.3,9)(13.3,9)(15.8,7.5)
\qbezier(15.8,7.5)(15.8,6)(15.8,6)
\qbezier(18.3,9)(18.3,9)(20.8,7.5)
\qbezier(15.8,6)(15.8,6)(18.3,4.5)

\end{picture}
\end{center}

\noindent
These values are for $T^{-d}_s(b)$.
The rigging for 
$b$ in question is defined to be their shift $+d\min(s,j)$
for length $j$ solitons,
leading to ($s=d=3$ in this example)

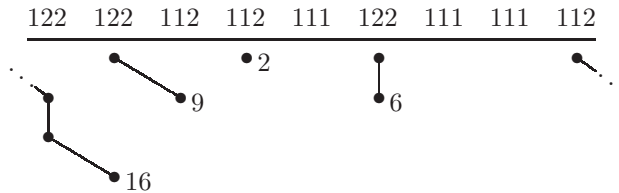
\begin{center}
\unitlength 10pt
\begin{picture}(24,8)(-0.2,4)
\put(0,10.2){$
122 \;\;\;\, 122 \;\;\;\, 112 \;\;\;\, 112 \;\;\;\, 111 \;\;\;\, 
122 \;\;\;\, 111 \;\;\;\, 111 \;\;\;\, 112$}
\put(0,9.7){\line(1,0){21.5}}

\put(3.3,9){\circle*{0.4}}
\put(8.3,9){\circle*{0.4}}\put(8.7,8.5){$2$}
\put(13.3,9){\circle*{0.4}}
\put(20.8,9){\circle*{0.4}}

\put(0.8,7.5){\circle*{0.4}}
\put(5.8,7.5){\circle*{0.4}}\put(6.2,7){$9$}
\put(13.3,7.5){\circle*{0.4}}\put(13.7,7){$6$}

\put(0.8,6){\circle*{0.4}}

\put(3.3,4.5){\circle*{0.4}}\put(3.7,4){$16$}

\qbezier(13.3,7.5)(13.3,7.5)(13.3,9)

\qbezier(5.8,7.5)(3.3,9)(3.3,9)

\qbezier(20.8,9)(20.8,9)(21.5,8.5)
\put(21.5,7.6){$\ddots$}

\qbezier(0.8,7.5)(0.8,7.5)(0.3,7.9)
\put(-0.8,7.8){$\ddots$}
\qbezier(0.8,6)(0.8,6)(0.8,7.5)
\qbezier(3.3,4.5)(3.3,4.5)(0.8,6)

\end{picture}
\end{center}

\noindent
Order the so obtained rigging for length $j$ solitons as
$r_{j,1}\le \cdots \le r_{j,m_j}$ and set 
\begin{equation}\label{eq:av2}
{\bf J} = (J_{j,\alpha})_{j \in H, 1 \le \alpha \le m_j},\quad 
\text{where }\;
J_{j,\alpha} =r_{j,\alpha}+\alpha-1.
\end{equation}
In the present example, $H=\{1,2,4\},\,
(m_1, m_2, m_4)=(1,2,1)$, $(p_1, p_2, p_4)=(1,4,9)$ and 
\begin{equation*}{\small
A=\begin{pmatrix}
p_1\!+\!m_1\!+\!1 & 2\min(1,2) & 2\min(1,2) & 2\min(1,4)\\ 
2\min(1,2) & p_2\!+\!m_2\!+\!3 & 3 & 2\min(2,4)\\ 
2\min(1,2) & 3 & p_2\!+\!m_2\!+\!3 & 2\min(2,4)\\ 
2\min(1,4) & 2\min(2,4) & 2\min(2,4) & p_4\!+\!m_4\!+\!7 
\end{pmatrix} = 
\begin{pmatrix}
3 & 2 & 2 & 2\\ 
2 & 9 & 3 & 4\\ 
2 & 3 & 9 & 4\\ 
2 & 4 & 4 & 17
\end{pmatrix},}
\end{equation*} so the angle variable is
\begin{equation}\label{eq:av1}
\begin{pmatrix}J_{1,1}\\ J_{2,1}\\ J_{2,2}\\ J_{4,1} \end{pmatrix}
=\begin{pmatrix}2\\ 6\\ 9+1 \\ 16 \end{pmatrix}\mod
A\Z^4 = \Z \begin{pmatrix}
3\\ 
2\\ 
2\\ 
2\end{pmatrix}
\oplus \Z
\begin{pmatrix}
2\\ 
9\\ 
3\\ 
4\end{pmatrix}
\oplus \Z
\begin{pmatrix}
2\\ 
3\\ 
9\\ 
4\end{pmatrix}
\oplus \Z
\begin{pmatrix}
2\\ 
4\\ 
4\\ 
17\end{pmatrix},
\end{equation}
where $+1$ is the contribution of $\alpha-1$
in $J_{j,\alpha} = r_{j,\alpha}+\alpha-1$.

This procedure for the direct scattering map $\Phi$
works provided that there is a cyclic shift $T_s^{-d}(b)$ 
such that no soliton stays across the boundary.
The case $s=1$ \cite{KTT} is such an example.
We conjecture it for general $s$ 
under the assumption (\ref{eq:pos}).

One can show that ${\bf J}\, \hbox{mod}\, \Gamma$ 
is independent of the possible non-uniqueness 
of such cyclic shifts.
The difference caused by such choices belong to $\Gamma$.
This can be observed, for example, by comparing 
(\ref{eq:av1}) and (\ref{eq:exx}).

Let us re-derive the result (\ref{eq:av1}) from (\ref{eq:Phi}).
The latter starts, for example, from the relation 
$b = T_3^{-5}(b_+)$, where 
\begin{equation*}
b_+ = 
\fbox{111} \otimes \fbox{122} \otimes 
\fbox{111} \otimes \fbox{111} \otimes \fbox{112}
\otimes \fbox{122} \otimes \fbox{122} \otimes \fbox{112} \otimes 
\fbox{112}
\end{equation*}
is a highest path corresponding to the rigged configuration
\begin{equation*}
\yngrc(4,{6},2,{3},2,{0},1,{1})
\end{equation*}
Thus (\ref{eq:Phi}) is evaluated as 
\begin{equation*}
{\bf I} - 5{\bf h}_3 = 
\begin{pmatrix}I_{1,1}\\ I_{2,1}\\ I_{2,2}\\ I_{4,1} \end{pmatrix}
=\begin{pmatrix}1\\ 0\\ 3+1 \\ 6 \end{pmatrix}-5
\begin{pmatrix}1\\ 2\\ 2 \\ 3 \end{pmatrix}
=\begin{pmatrix}-4\\ -10\\ -6 \\ -9 \end{pmatrix}
\mod A\Z^4.
\end{equation*}
This certainly coincides with the result (\ref{eq:av1}) since 
the difference 
\begin{equation}\label{eq:exx}
\begin{pmatrix}2\\ 6\\ 10 \\ 16 \end{pmatrix} - 
\begin{pmatrix}-4\\ -10\\ -6 \\ -9 \end{pmatrix} = 
\begin{pmatrix}6\\ 16\\ 16 \\ 25 \end{pmatrix}
=\begin{pmatrix}
2\\ 
9\\ 
3\\ 
4\end{pmatrix}
+
\begin{pmatrix}
2\\ 
3\\ 
9\\ 
4\end{pmatrix}
+
\begin{pmatrix}
2\\ 
4\\ 
4\\ 
17\end{pmatrix}
\end{equation}
belongs to $A\Z^4$.

\subsection{Inverse scattering}\label{ssec:ism}

According to Theorem \ref{th:cd}, 
the dynamics of the generalized periodic box-ball system 
is transformed to a straight motion 
(\ref{eq:tla}) in the set ${\mathcal J}(m)$ 
of angle variables.
To complete the inverse scattering method, 
one needs the inverse scattering map $\Phi^{-1}$ 
from ${\mathcal J}(m)$ back to paths $\Pth(m)$.
Under the condition (\ref{eq:pos}), it is easy to 
show that any element 
${\bf I} \in {\mathcal J}(m)$ 
has a (not necessarily unique) representative form 
${\bf I} = \sum_ld_l{\bf h}_l+
(r_{j,\alpha}+\alpha-1)_{j \in H,1 \le \alpha \le m_j}$
such that $0 \le r_{j,1} \le \cdots \le r_{j,m_j} \le p_j$
by using the equivalence under $\Gamma$. 
If $\mu$ denotes the Young diagram for $m=\{m_k\}$
and $r = (r_{j,\alpha})_{j \in H, 1 \le \alpha \le m_j}$,
then $(\mu,r)$ becomes a rigged configuration.
Letting $b_+$ be the highest path corresponding to it,
$\Phi^{-1}({\bf I}) := \bigl(\prod_lT_l^{d_l}\bigr)(b_+)$ 
is independent of the choice of the representative form and 
yields the inverse of $\Phi$.
Actually, one can always take $d_l=0$ for $l \neq 1$.

Our solution of the initial value problem is achieved by
the commutative diagram (\ref{eq:cd}), namely 
the composition:
\begin{equation*}
\Pth(m) \overset{\Phi}{\longrightarrow} {\mathcal J}(m) 
\overset{\{T_l\} }{\longrightarrow} {\mathcal J}(m) 
\overset{\Phi^{-1}}{\longrightarrow} \Pth(m),
\end{equation*}
where the number of computational steps is 
independent of the time evolution.
As an illustration we derive 
\begin{align*}
T_2^{1000}(b)
&=\fbox{111} \otimes \fbox{111} \otimes \fbox{112} \otimes 
\fbox{112} \otimes \fbox{122} \otimes \fbox{122} \otimes 
\fbox{112} \otimes \fbox{111} \otimes \fbox{122},\\
T_4^{1000}(b) 
&= \fbox{112} \otimes \fbox{112} \otimes \fbox{112} \otimes 
\fbox{111} \otimes \fbox{122} \otimes \fbox{111} \otimes 
\fbox{111} \otimes \fbox{122} \otimes \fbox{122}
\end{align*}
for $b$ given in (\ref{eq:bex}).
{}From (\ref{eq:av1}) and (\ref{eq:tla}),
the angle variables for $T^{1000}_2(b)$ and $T_4^{1000}(b)$ 
are given by 
\begin{align*}
\begin{pmatrix}1002\\ 2006\\ 2010 \\ 2016 \end{pmatrix}
&= \begin{pmatrix}0+82\\ 0+82\\ 3+1+82 \\ 6+82 \end{pmatrix}
+108 \begin{pmatrix}
3\\ 
2\\ 
2\\ 
2\end{pmatrix}
+ 129
\begin{pmatrix}
2\\ 
9\\ 
3\\ 
4\end{pmatrix}
+ 129
\begin{pmatrix}
2\\ 
3\\ 
9\\ 
4\end{pmatrix}
+ 40
\begin{pmatrix}
2\\ 
4\\ 
4\\ 
17\end{pmatrix},\\
\begin{pmatrix}1002\\ 2006\\ 2010 \\ 4016 \end{pmatrix}
&= \begin{pmatrix}0+222\\ 0+222\\ 3+1+222 \\ 6+222 \end{pmatrix}
+28 \begin{pmatrix}
3\\ 
2\\ 
2\\ 
2\end{pmatrix}
+ 84
\begin{pmatrix}
2\\ 
9\\ 
3\\ 
4\end{pmatrix}
+ 84
\begin{pmatrix}
2\\ 
3\\ 
9\\ 
4\end{pmatrix}
+ 180
\begin{pmatrix}
2\\ 
4\\ 
4\\ 
17\end{pmatrix}.
\end{align*}
In the right hand sides, the last four terms belong to $A\Z^4$, 
hence can be neglected.
The first terms correspond to $T_1^{82}$ and $T_1^{222}$ of the 
rigged configuration ($+1$ is removed as the ``$\alpha-1$ part")
\begin{equation*}
\yngrc(4,{6},2,{3},2,{0},1,{0}) 
\end{equation*}
which is mapped, under the KKR bijection, to the highest path
\begin{equation*}
b'_+:=\fbox{111} \otimes \fbox{122} \otimes \fbox{111} \otimes 
\fbox{112} \otimes \fbox{111} \otimes \fbox{122} \otimes 
\fbox{122} \otimes \fbox{112} \otimes \fbox{112}. 
\end{equation*}
One can check $T^{82}_1(b'_+) = T^{1000}_2(b)$
and $T^{222}_1(b'_+) = T^{1000}_4(b)$ completing the 
derivation.

\subsection{State counting and periodicity}\label{ssec:pb}
The matrix $A$ (\ref{eq:a}) was originally introduced in the study 
of the Bethe equation at $q=0$ \cite{KN}.
{}From this connection we have
\begin{theorem}[\cite{KN} Theorems 3.5, 4.9]\label{th:kn}
\begin{align}
\vert {\mathcal J}(m) \vert &= 
(\det F)\prod_{j \in H} \frac{1}{m_j}
\binom{p_j + m_j - 1}{m_j - 1}.\label{eq:omega}\\
F &= (F_{j,k})_{j,k \in H},\quad 
F_{j,k} = \delta_{j,k}p_j + 2\min(j,k)m_k.\label{eq:f}
\end{align}
\end{theorem}
Combined with Theorem \ref{th:cd}, this yields a formula for 
$\vert \Pth(m) \vert$, namely, the number of states
characterized by the conserved quantity.
For $B_3^{\otimes 9}$ and 
$m$ corresponding to the Young diagram $(4,2,2,1)$ 
considered above, we have  $\vert {\mathcal J}(m) \vert=990$.

{}From (\ref{eq:av}), all the paths $b \in \Pth(m)$ obey 
the relation 
\begin{equation*}
T_l^{{\mathcal N}_l}(b) = b\quad
\hbox{if}\;\;
{\mathcal N}_l{\bf h}_l \in \Gamma.
\end{equation*}
Writing ${\mathcal N}_l{\bf h}_l = A{\bf n}$, components of
the vector ${\bf n}=(n_{j,\alpha})$ are given by
\begin{equation*}
n_{j,\alpha} = {\mathcal N}_l \frac{\det A[j\alpha]}{\det A},
\end{equation*}
where  $A[j\alpha]$ is obtained by replacing $(j\alpha)$th 
column of $A$ by ${\bf h}_l$.
It is elementary to check
\begin{equation*}
\frac{\det A[j\alpha]}{\det A} = \frac{\det F[j]}{\det F}.
\end{equation*}
$F[j]$ is obtained from $F$ by replacing 
the $j$th column by the $l$-dependent vector $(\min(l,k))_{k \in H}$.
The independence on $\alpha$
reflects the symmetry of $A$ (\ref{eq:a}) within blocks.
Thus the {\em generic} period of $\Pth(m)$, 
namely the minimum ${\mathcal N}_l$ such that  
${\mathcal N}_l{\bf h}_l \in \Gamma$ is
\begin{equation}\label{eq:lcm}
{\mathcal N}_l = {\rm LCM}\!\left(
\frac{\det F}{\det F[j_1]}, \ldots, \frac{\det F}{\det F[j_g]}
\right),
\end{equation}
where by ${\rm LCM}(r_1, \ldots, r_g)$ for 
rational numbers $r_1, \ldots, r_g$, we mean  
the smallest positive integer  
in $\Z r_1 \cap \cdots \cap \Z r_g$.
When $\det F[j_k]=0$, the entry 
${\det F}/{\det F[j_k]}$ is to be excluded. 
For $B^{\otimes 9}_3$ and the Young diagram $(4,2,2,1)$ 
in the above example, we have
\begin{equation*}
F=\begin{pmatrix}
3 & 4 & 2 \\
2 & 12 & 4 \\
2 & 8 & 17
\end{pmatrix},\;
F[1]=\begin{pmatrix}
1 & 4 & 2 \\
1 & 12 & 4 \\
1 & 8 & 17
\end{pmatrix},\;
F[2]=\begin{pmatrix}
3 & 1 & 2 \\
2 & 1 & 4 \\
2 & 1 & 17
\end{pmatrix},\;
F[4]=\begin{pmatrix}
3 & 4 & 1 \\
2 & 12 & 1 \\
2 & 8 & 1
\end{pmatrix}
\end{equation*}
for $l=1$, leading to 
${\mathcal N}_1 = 
{\rm LCM}\left(\frac{99}{28}, \frac{396}{13}, 99\right)=396$. 
Similar calculations yield
\begin{equation*}
{\mathcal N}_1=396,\quad
{\mathcal N}_2=99,\quad
{\mathcal N}_3=9,\quad
{\mathcal N}_l=11\;(l \ge 4).
\end{equation*}
$T_l^{{\mathcal N}_l}(b)=b$ can be directly checked
for $b$ in (\ref{eq:bex}).
In fact $T_4^{11}(b)=b$ has been demonstrated 
in the introduction.
For the fundamental period,
formally the same closed formula 
as eq.(4.26) in \cite{KTT} is valid.
The formula (\ref{eq:lcm}) agrees with
the most general conjecture on $A^{(1)}_n$ \cite{KT}.
For $s=1$ and $l=\infty$, it was originally obtained 
in \cite{YYT} by a combinatorial argument. 

\subsection{Ultradiscrete Riemann theta lattice}\label{ssec:udr}
Let us present an explicit formula for the inverse 
scattering map 
$\Phi^{-1}: {\mathcal J}(m) \rightarrow \Pth(m)$
in terms of the ultradiscrete Riemann theta function.
We keep assuming the condition (\ref{eq:pos}).

For general $\{m_j\}$, Theorem 5.1 in \cite{KS2}
remains valid if the vacancy number $p_j$ is
replaced by (\ref{eq:va}) in this paper.
Here we restrict ourselves to the case 
$m_j=1$ for all $j \in H = \{j_1 < \cdots < j_g\}$ for simplicity.
Thus ${\mathcal J}(m)$ (\ref{eq:av}) 
and $A$ (\ref{eq:a}) reduce to
\begin{align}
{\mathcal J}(m) &= \Z^g/A\Z^g,\label{eq:js}\\
A&=(A_{j,k})_{j,k \in H},\quad
A_{j,k} = \delta_{j,k}p_j + 2\min(j,k).\nonumber
\end{align}
Following \cite{KS1}, we introduce the ultradiscrete 
Riemann theta function by
\begin{equation}\label{eq:urt}
\begin{split}
\Theta({\bf z}) &= \lim_{\epsilon \rightarrow +0}
\epsilon\log\left( \,\sum_{{\bf n} \in \Z^g}
\exp\Bigl(-\frac{{}^t{\bf n}A{\bf n}/2+{}^t{\bf n}{\bf z}}{\epsilon}
\Bigr)\right)\\
&= -\min_{{\bf n} \in \Z^g}
\{{}^t{\bf n}A{\bf n}/2+{}^t{\bf n}{\bf z}\},
\end{split}
\end{equation}
which enjoys the quasi-periodicity:
\begin{equation}\label{eq:qp}
\Theta({\bf z} + {\bf v}) = 
{}^t{\bf v}A^{-1}({\bf z} + {\bf v}/2) + \Theta({\bf z})\quad 
\hbox{for any } \; {\bf v} \in 
\Gamma = A\Z^g.
\end{equation}

Consider the planar square lattice 
\begin{equation*}
\begin{picture}(180,140)(25,35)
\setlength{\unitlength}{1mm}

\put(2,45){\line(1,0){51}}\put(56,44.8){$\ldots$}\put(63,45){\line(1,0){17}}
\put(2,30){\line(1,0){51}}\put(56,29.8){$\ldots$}\put(63,30){\line(1,0){17}}

\put(13,58){$B_{s_1}$}\put(28, 58){$B_{s_2}$}\put(43, 58){$B_{s_3}$}
\put(68, 58){$B_{s_L}$}

\put(-4,44){$B_{l_1}$}\put(-4,29){$B_{l_2}$}

\put(15,20){\line(0,1){35}}\put(30,20){\line(0,1){35}}\put(45,20){\line(0,1){35}}
\put(70,20){\line(0,1){35}}

\put(6,50){${\bf z}_{0,0}$}\put(20,50){${\bf z}_{0,1}$}\put(35,50){${\bf z}_{0,2}$}
\put(74,50){${\bf z}_{0,L}$}

\put(6,36.5){${\bf z}_{1,0}$}\put(20,36.5){${\bf z}_{1,1}$}
\put(35,36.5){${\bf z}_{1,2}$}
\put(74,36.5){${\bf z}_{1,L}$}

\put(6,22.5){${\bf z}_{2,0}$}
\put(20,22.5){${\bf z}_{2,1}$}
\put(35,22.5){${\bf z}_{2,2}$}
\put(74,22.5){${\bf z}_{2,L}$}

\put(14.5,15){$\vdots$}\put(29.3,15){$\vdots$}\put(44.5,15){$\vdots$}
\put(69.5,15){$\vdots$}

\end{picture}
\end{equation*}
where we assume $s_1=\cdots = s_L = s$ for the time being.
${\bf z}_{t,k}$ is defined by
\begin{equation}
{\bf z}_{t,k} = {\bf I} - \frac{\bf p}{2} 
+  {\bf h}_{l_1} + \cdots +  
{\bf h}_{l_t }- {\bf h}_{s_1} - \cdots - {\bf h}_{s_k},
\end{equation}
where ${\bf I} \in {\mathcal J}(m)$ is the angle variable 
of a path 
$b = b_1 \otimes \cdots \otimes b_L 
\in B_{s_1} \otimes \cdots \otimes B_{s_L}$, and 
${\bf p}={}^t(p_{j_1},\ldots, p_{j_g})$.
To each edge, either vertical or horizontal, 
we attach a number $\Theta({\bf z}, {\bf z}')$ via the rule:
\begin{equation}\label{eq:zt}
\begin{picture}(140,40)(50,42)
\setlength{\unitlength}{1mm}

\put(-14.5,19){${\bf z}'$}
\put(-10,15){\line(0,1){9}}
\put(-7.5,19){${\bf z}$\;\;, }

\put(6,22.5){${\bf z}$}
\put(1.5,20){\line(1,0){10}}
\put(6,15.5){${\bf z}'$}
\put(16,19){$\rightarrow \;\;
\Theta({\bf z}, {\bf z}') := \Theta({\bf z}) - \Theta({\bf z}') - 
\Theta({\bf z}+{\bf h}_\infty) + \Theta({\bf z}'+{\bf h}_\infty)$.}
\end{picture}
\end{equation}
The rule (\ref{eq:zt}) can also be described neatly 
in terms of a two-layer cubic lattice whose sites 
are assigned with ${\bf z}_{t,k}$ or ${\bf z}_{t,k} + {\bf h}_\infty$.
Assign the 2-component integer vectors 
\begin{equation}\label{eq:xyt}
\begin{split}
x_{t,k}&=\Bigl(l_t - \Theta({\bf z}_{t-1,k-1}, {\bf z}_{t, k-1}),
\Theta({\bf z}_{t-1,k-1}, {\bf z}_{t, k-1})\Bigr),\\
y_{t,k}&=\Bigl(s_k-\Theta({\bf z}_{t-1,k}, {\bf z}_{t-1,k-1}), 
\Theta({\bf z}_{t-1,k}, {\bf z}_{t-1,k-1})\Bigr)
\end{split}
\end{equation}
to the edges as follows:
\begin{equation*}
\begin{picture}(60,85)(10,-10)
\setlength{\unitlength}{1mm}

\put(8,23){$y_{t,k}$}
\put(-11,9){$x_{t,k}$}
\put(26,9){$x_{t,k+1}$}
\put(7.5,-3){$y_{t+1,k}$}

\put(-4,10){\line(1,0){28}}
\put(10,0){\line(0,1){20}}

\put(-3.5,14.5){${\bf z}_{t-1,k-1}$} \put(13.5,14.5){${\bf z}_{t-1,k}$}
\put(-1,4){${\bf z}_{t,k-1}$} \put(13.5,4){${\bf z}_{t,k}$}

\end{picture}
\end{equation*}
By the same argument as \cite{KS1}, one can show that 
the ultradiscrete tau function in \cite{KSY} 
for the periodic system coincides essentially 
with $\Theta$ (\ref{eq:urt}) here.
{}From this fact and Theorem \ref{th:cd} we obtain
\begin{theorem}\label{th:xy}
\begin{equation}\label{eq:rxy}
x_{t,k} \in B_{l_t},\; y_{t,k} \in B_{s_k},\quad 
R(x_{t, k} \otimes y_{t, k}) = y_{t+1, k} \otimes x_{t, k+1},
\end{equation}
where $R$ denotes the combinatorial $R$.
The path $b$ is reproduced from 
${\bf I} \in {\mathcal J}(m)$ by
\begin{equation*}
b = y_{1,1} \otimes y_{1,2} \otimes \cdots \otimes y_{1,L}.
\end{equation*}
\end{theorem}
The assertion $x_{t,k} \in B_{l_t}$ means that
$\Theta({\bf z}_{t-1,k-1}, {\bf z}_{t, k-1})$ is an integer 
in the range $[0, l_t]$, and similarly for $y_{t,k} \in B_{s_k}$.
The periodic boundary condition in the horizontal direction 
\begin{equation*}
x_{t,0} = x_{t, L}
\end{equation*}
is valid for any $t$.
This can easily be checked by using 
${\bf z}_{t,0}-{\bf z}_{t,L} = {\bf h}_{s_1} + \cdots 
+ {\bf h}_{s_L} = A{}^t(1,1\ldots,1)$ and the quasi-periodicity 
(\ref{eq:qp}). 
Theorem \ref{th:xy} tells that 
$x_{t,k}$ and $y_{t,k}$ obey the local dynamics (combinatorial $R$) 
of the generalized periodic box-ball system.
As a result, we get
\begin{equation}\label{eq:sol}
T_{l_t}\cdots T_{l_1}(b) 
= y_{t+1,1} \otimes y_{t+1,2} \otimes \cdots \otimes y_{t+1,L}.
\end{equation}
In this way the solution of the initial value problem 
under arbitrary time evolutions $\{T_{l_t}\}$
is written down explicitly.
Note that we have given an explicit 
formula not only for $y_{t,k}$ but also $x_{t,k}$
which is often called  ``carrier" \cite{TM}.

These aspects of the periodic box-ball system on  
$\Pth = B_s^{\otimes L}$ 
admit a generalization to the inhomogeneous case
$\Pth = B_{s_1} \otimes \cdots \otimes B_{s_L}$.
A typical example is a combinatorial version 
of Yang's system \cite{Ya,V}, which is endowed with 
the family of time evolutions 
$\{T_{s_1}, \ldots, T_{s_L}\}$.
The well-known relation 
$T_{s_1}T_{s_2}\cdots T_{s_L} = {\rm Id}$ is understood from 
${\bf h}_{s_1} + \cdots 
+ {\bf h}_{s_L} = A{}^t(1,1\ldots,1) \in \Gamma$
in our linearization scheme.

To adapt the formalism to the inhomogeneous situation 
$\Pth = B_{s_1} \otimes \cdots \otimes B_{s_L}$, 
we employ (\ref{eq:av2}) to specify $\Phi$
and replace the vacancy number (\ref{eq:va}) with 
\begin{equation*}
p_j = \sum_{i=1}^L\min(s_i,j) - 2E_j.
\end{equation*}
Under the condition (\ref{eq:pos}), we conjecture 
that Theorem \ref{th:xy}  remains valid.
For an illustration, we consider the path 
\begin{equation}\label{eq:pex}
b=11 \otimes 2 \otimes 1 \otimes 2 \otimes 122 \otimes
1 \otimes 12 \otimes 2 \otimes 1.
\end{equation}
The action-angle variable is depicted as
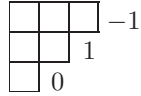
\begin{equation*}
\begin{picture}(50,50)(30,0)
\put(0,22){$\yng(3,2,1)$}
\put(37,31){$-1$}
\put(28,20){$1$}
\put(16,8){$0$}
\end{picture}
\end{equation*}
So the vacancy numbers ${\bf p}$,
the period matrix $A$ and the angle variable
${\bf I}$ are given by 
\begin{equation*}
{\bf p} = \begin{pmatrix}
p_1 \\ p_2 \\ p_3 
\end{pmatrix}
= \begin{pmatrix}
3 \\ 2 \\ 1 
\end{pmatrix},\quad
A= \begin{pmatrix}
5 & 2 & 2 \\ 2 & 6 & 4 \\ 2 & 4 & 7 
\end{pmatrix},\quad
{\bf I} = \begin{pmatrix}
0 \\ 1 \\ -1 
\end{pmatrix}.
\end{equation*}
We set $(s_1,\ldots, s_9) = (2,1,1,1,3,1,2,1,1)$ according to (\ref{eq:pex}).
Let us take $(l_1, l_2, l_3)=(2,1,3)$ and consider the 
time evolution 
$b \rightarrow T_{l_1}(b) \rightarrow  T_{l_2}T_{l_1}(b) 
\rightarrow T_{l_3}T_{l_2}T_{l_1}(b)$.
Then the edge variables 
exhibit the following pattern:
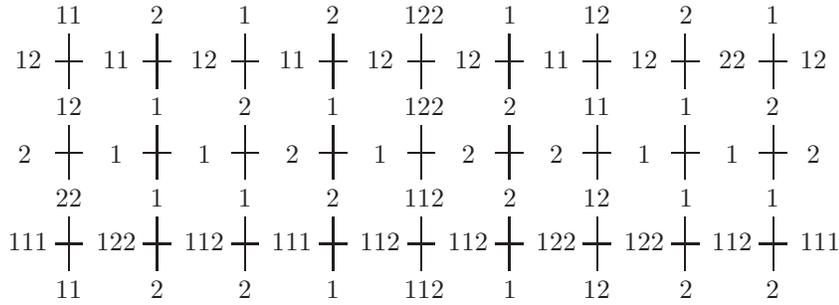
\begin{equation*}
\begin{picture}(200,140)(65,-30)
\setlength{\unitlength}{0.9mm}

\multiput(0,0)(0,13.5){3}{
\multiput(7,0.5)(13,0){9}{
\put(0,2.5){\line(1,0){4}}\put(2,-1.5){\line(0,1){8}}
}}

\put(7,35.5){11}\put(21,35.5){2}\put(34,35.5){1}\put(47,35.5){2}
\put(58.5,35.5){122}\put(73.2,35.5){1}\put(85,35.5){12}\put(99.2,35.5){2}
\put(112,35.5){1}

\put(7,22){12}\put(21,22){1}\put(34,22){2}\put(47,22){1}
\put(58.5,22){122}\put(73.2,22){2}\put(85,22){11}\put(99.2,22){1}
\put(112,22){2}

\put(7,8.5){22}\put(21,8.5){1}\put(34,8.5){1}\put(47,8.5){2}
\put(58.5,8.5){112}\put(73.2,8.5){2}\put(85,8.5){12}\put(99.2,8.5){1}
\put(112,8.5){1}

\put(7,-5){11}\put(21,-5){2}\put(34,-5){2}\put(47,-5){1}
\put(58.5,-5){112}\put(73.2,-5){1}\put(85,-5){12}\put(99.2,-5){2}
\put(112,-5){2}

\put(1,29){12}\put(14,29){11}\put(27,29){12}\put(40,29){11}\put(53,29){12}
\put(66,29){12}\put(79,29){11}\put(92,29){12}\put(105,29){22}\put(117,29){12}

\put(1.5,15){2}\put(15,15){1}\put(28,15){1}\put(41,15){2}\put(54,15){1}
\put(67,15){2}\put(80,15){2}\put(93,15){1}\put(106,15){1}\put(118,15){2}

\put(0,2){111}\put(13,2){122}\put(26,2){112}
\put(39,2){111}\put(52,2){112}\put(65,2){112}\put(78,2){122}
\put(91,2){122}\put(104,2){112}\put(117,2){111}

\end{picture}
\end{equation*}
In terms of the edge variable $\Theta({\bf z}, {\bf z}')$
(\ref{eq:zt}) 
representing the number of the letter 2 in tableaux on edges,
this looks as
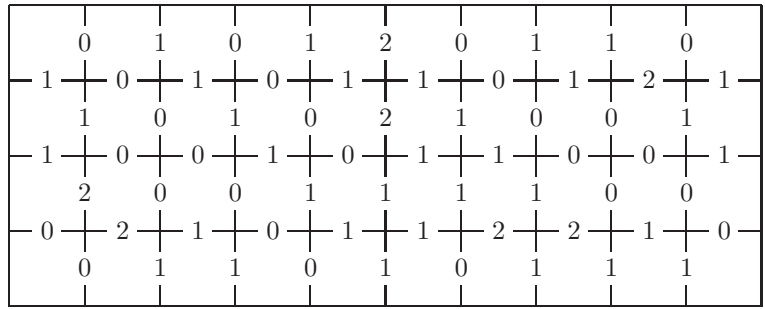
\begin{equation}\label{tab:ev}
\begin{picture}(120,130)(90,-10)
\setlength{\unitlength}{1mm}
\multiput(-0.2,0)(0,40){2}{
\put(0,0){\line(1,0){100}}
}
\multiput(-0.2,0)(100,0){2}{
\put(0,0){\line(0,1){40}}
}

\multiput(-0.2,0)(0,10){3}{
\multiput(0,0)(10,0){10}{
\put(0,10){\line(1,0){3}}\put(10,10){\line(-1,0){3}}
}}

\multiput(-0.2,0)(0,10){4}{
\multiput(0,0)(10,0){10}{
\put(10,0){\line(0,1){2.7}}\put(10,10){\line(0,-1){2.7}}
}}

\put(9,34){0}\put(19,34){1}\put(29,34){0}\put(39,34){1}
\put(49,34){2}\put(59,34){0}\put(69,34){1}\put(79,34){1}\put(89,34){0}

\put(9,24){1}\put(19,24){0}\put(29,24){1}\put(39,24){0}
\put(49,24){2}\put(59,24){1}\put(69,24){0}\put(79,24){0}\put(89,24){1}

\put(9,14){2}\put(19,14){0}\put(29,14){0}\put(39,14){1}
\put(49,14){1}\put(59,14){1}\put(69,14){1}\put(79,14){0}\put(89,14){0}

\put(9,4){0}\put(19,4){1}\put(29,4){1}\put(39,4){0}
\put(49,4){1}\put(59,4){0}\put(69,4){1}\put(79,4){1}\put(89,4){1}

\put(4,29){1}\put(14,29){0}\put(24,29){1}\put(34,29){0}\put(44,29){1}
\put(54,29){1}\put(64,29){0}\put(74,29){1}\put(84,29){2}\put(94,29){1}

\put(4,19){1}\put(14,19){0}\put(24,19){0}\put(34,19){1}\put(44,19){0}
\put(54,19){1}\put(64,19){1}\put(74,19){0}\put(84,19){0}\put(94,19){1}

\put(4,9){0}\put(14,9){2}\put(24,9){1}\put(34,9){0}\put(44,9){1}
\put(54,9){1}\put(64,9){2}\put(74,9){2}\put(84,9){1}\put(94,9){0}

\end{picture}
\end{equation}
The table of the values of the ultradiscrete 
Riemann theta $\Theta({\bf z}_{t,k})$ is as follows:
\begin{equation*}
\begin{picture}(120,130)(90,-10)
\setlength{\unitlength}{1mm}
\multiput(0,0)(0,10){5}{
\put(0,0){\line(1,0){100}}
}
\multiput(0,0)(10,0){11}{
\put(0,0){\line(0,1){40}}
}

\put(4,34){0}\put(14,34){0}\put(24,34){1}\put(34,34){2}\put(44,34){4}
\put(54,34){8}\put(63,34){10}\put(73,34){14}\put(83,34){17}\put(93,34){20}

\put(4,24){0}\put(14,24){0}\put(24,24){0}\put(34,24){1}\put(44,24){2}
\put(54,24){5}\put(64,24){7}\put(73,24){10}\put(83,24){12}\put(93,24){15}

\put(4,14){0}\put(14,14){0}\put(24,14){0}\put(34,14){0}\put(44,14){1}
\put(54,14){3}\put(64,14){5}\put(74,14){8}\put(83,14){10}\put(93,14){12}

\put(4,4){2}\put(14,4){0}\put(24,4){0}\put(34,4){0}\put(44,4){0}
\put(54,4){1}\put(64,4){2}\put(74,4){4}\put(84,4){6}\put(94,4){8}

\end{picture}
\end{equation*}
Similarly the values of $\Theta({\bf z}_{t,k}+{\bf h}_\infty)$ are given as follows:
\begin{equation*}
\begin{picture}(120,130)(90,-10)
\setlength{\unitlength}{1mm}
\multiput(0,0)(0,10){5}{
\put(0,0){\line(1,0){100}}
}
\multiput(0,0)(10,0){11}{
\put(0,0){\line(0,1){40}}
}

\put(4,34){0}\put(14,34){0}\put(24,34){0}\put(34,34){1}\put(44,34){2}
\put(54,34){4}\put(64,34){6}\put(74,34){9}\put(83,34){11}\put(93,34){14}

\put(4,24){1}\put(14,24){0}\put(24,24){0}\put(34,24){0}\put(44,24){1}
\put(54,24){2}\put(64,24){3}\put(74,24){6}\put(84,24){8}\put(93,24){10}

\put(4,14){2}\put(14,14){0}\put(24,14){0}\put(34,14){0}\put(44,14){0}
\put(54,14){1}\put(64,14){2}\put(74,14){4}\put(84,14){6}\put(94,14){8}

\put(4,4){4}\put(14,4){2}\put(24,4){1}\put(34,4){0}\put(44,4){0}
\put(54,4){0}\put(64,4){1}\put(74,4){2}\put(84,4){3}\put(94,4){4}
\end{picture}
\end{equation*}
{}From these values of 
$\Theta({\bf z}_{t,k})$ and $\Theta({\bf z}_{t,k}+{\bf h}_\infty)$,
the table (\ref{tab:ev}) is reproduced by the rule (\ref{eq:zt}). 

\section{Summary}\label{sec:sum}
We have introduced the local energy distribution 
of paths in section \ref{sec:led} and reformulated the 
KKR map $\phi$ in Theorem \ref{s:main}.
Combined with the result for $\phi^{-1}$ \cite{KOSTY,Sa},
it completes the crystal interpretation of the 
KKR bijection for $U_q(\widehat{\mathfrak{sl}}_2)$.

In section \ref{sec:pbbs} the generalized periodic 
box-ball system on $B^{\otimes L}_s$ is studied.
Under the condition (\ref{eq:pos}),  the set of paths
${\widehat {\mathcal P}}(m)$ (\ref{eq:pm}) 
characterized by conserved quantities  
enjoy all the properties (i) and (ii) 
stated under (\ref{eq:inv}).
The action-angle variables are introduced in section 
\ref{ssec:aav}.
The inverse scattering formalism, i.e., 
simultaneous linearization of the commuting family of 
time evolutions, is established in 
Theorem \ref{th:cd}.
It leads to the formulas for state counting (Theorem \ref{th:kn}),
generic period (\ref{eq:lcm}), and an algorithm for solving 
the initial value problem.
According to Theorem \ref{th:xy},
the solution of the initial value problem (\ref{eq:sol}) is expressed 
in terms of the ultradiscrete Riemann theta function (\ref{eq:urt}).
A similar formula has been obtained also for the carrier variable 
$x_{t,k}$ simultaneously.
These results extend the $s=1$ case \cite{KTT,KS1} and 
agree with the conjecture on the most general case \cite{KT}.

\appendix
\section{Crystals and Combinatorial $R$}\label{sec:cry}
\subsection{Crystals}
We recapitulate the basic facts in the crystal base theory 
\cite{K,KMN,NY}.
Let $B_l$ the the crystal of the $l$-fold symmetric 
tensor representation of $U_q(A^{(1)}_1)$.
As a set it is given by 
$B_l = \{x=(x_1,x_2) \in (\Z_{\ge 0})^2 \mid x_1+x_2=l\}$.
The element $(x_1,x_2)$ will also be expressed as 
the length $l$ row shape semistandard tableau 
containing the letter $i$ $x_i$ times.
For example,
$B_1=\big\{\fbox{1},\fbox{2}\big\}, 
B_2 = \big\{\fbox{11}, \fbox{12}, \fbox{22} \big\}$.
(We shall often omit the frames of the tableaux.)
The action of Kashiwara operators 
${\tilde f}_i, {\tilde e}_i: B \rightarrow B \sqcup \{0\}\, (i=0,1)$ reads 
$({\tilde f}_ix)_j = x_j-\delta_{j,i}+\delta_{j,i+1}$ and 
$({\tilde e}_ix)_j = x_j+\delta_{j,i}-\delta_{j,i+1}$,
where all the indices are in $\Z_2$, and if the result does not 
belong to $(\Z_{\ge 0})^2$, it should be understood as $0$.
The classical part of the weight of $x=(x_1,x_2) \in B_l$ is 
$\wt(x)=l\Lambda_1-x_2\alpha_1 = (x_1-x_2)\Lambda_1$,
where $\Lambda_1$ and $\alpha_1=2\Lambda_1$ are 
the fundamental weight and the simple root of $A_1$.

For any $b\in B$, set 
\[
\veps_i(b)=\max\{m\ge0\mid \et{i}^m b\ne0\},\quad
\vphi_i(b)=\max\{m\ge0\mid \ft{i}^m b\ne0\}.
\]
By the definition one has $\veps_i(x)=x_{i+1}$ and $\vphi_i(x)=x_i$
for $x=(x_1,x_2) \in B_l$.
Thus $\veps_i(x)+\vphi_i(x)=l$ is valid for any $x \in B_l$.

For two crystals $B$ and $B'$, one can define the tensor product
$B\ot B'=\{b\ot b'\mid b\in B,b'\in B'\}$. 
The operators $\et{i},\ft{i}$ act 
on $B\ot B'$ by
\begin{eqnarray}
\et{i}(b\ot b')&=&\left\{
\begin{array}{ll}
\et{i} b\ot b'&\mbox{ if }\vphi_i(b)\ge\veps_i(b')\\
b\ot \et{i} b'&\mbox{ if }\vphi_i(b) < \veps_i(b'),
\end{array}\right. \label{eq:erule}\\
\ft{i}(b\ot b')&=&\left\{
\begin{array}{ll}
\ft{i} b\ot b'&\mbox{ if }\vphi_i(b) > \veps_i(b')\\
b\ot \ft{i} b'&\mbox{ if }\vphi_i(b)\le\veps_i(b').
\end{array}\right. \label{eq:frule}
\end{eqnarray}
Here $0\ot b'$ and $b\ot 0$ should be understood as $0$. 
The tensor product $B_{l_1}\otimes \cdots \otimes B_{l_k}$
is obtained by repeating the above rule.
The classical part of the weight of $b\in B$ for any 
$B=B_{l_1}\otimes \cdots \otimes B_{l_k}$ is given by 
$\wt(b) = (\varphi_1(b)-\varepsilon_1(b))\Lambda_1
=(\varepsilon_0(b)-\varphi_0(b))\Lambda_1$.

Let $s_i\, (i=0,1)$ be the Weyl group operator \cite{K}
acting on any crystal $B$ as
\begin{equation}\label{eq:wdef}
s_i(b) = \begin{cases}
{\tilde f}_i^{\varphi_i(b)-\varepsilon_i(b)}(b) & 
\varphi_i(b) \ge \varepsilon_i(b),\\
{\tilde e}_i^{\varepsilon_i(b)-\varphi_i(b)}(b) & 
\varphi_i(b) \le \varepsilon_i(b)
\end{cases}
\end{equation}
for $b \in B$.
Let 
\begin{equation}\label{eq:odef}
\omega: (x_1,x_2) \mapsto (x_2,x_1)
\end{equation}
be the involutive Dynkin diagram automorphism of $B_l$.
We extend it to any $B = B_{l_1}\otimes \cdots \otimes B_{l_k}$ 
by $\omega(B) 
= \omega(B_{l_1})\otimes \cdots \otimes \omega(B_{l_k})$.
Then ${\widehat W}(A^{(1)}_1)
= \langle \omega, s_0, s_1\rangle$ acts on
$B$ as the extended affine Weyl group of type $A^{(1)}_1$.

The action of ${\tilde f}_i, {\tilde e}_i$ and $s_i$ 
is determined in principle by (\ref{eq:erule}) and (\ref{eq:frule}).
Here we explain the {\em signature rule} to find the action 
on any $B_{l_1} \otimes \cdots \otimes B_{l_k}$, 
which is of great practical use.
The $i$-signature of an element $b \in B_l$ is the symbol
$\overbrace{-\cdots -}^{\varepsilon_i(b)}
\overbrace{+ \cdots +}^{\varphi_i(b)}$.
The $i$-signature of 
the tensor product $b_1 \otimes \cdots \otimes b_k \in 
B_{l_1} \otimes \cdots \otimes B_{l_k}$ is the array of 
the $i$-signature of each $b_j$.
Here is an example from 
$B_5\otimes B_2\otimes B_1 \otimes B_4$:
\begin{equation*}
\underset{1-{\rm signature}}{{\underset{0-{\rm signature}}
{\phantom{1-signature}}}}
\underset{-++++}{\underset{----+}{11112}} \otimes 
\underset{-+}{\underset{-+}{12}} \otimes 
\underset{-}{\underset{+}{2}} \otimes 
\underset{--++}{\underset{--++}{1122}}
\end{equation*}
where $1122$ for example represents $\fbox{1122} \in B_4$ and 
not $\fbox{1}\otimes \fbox{1} \otimes 
\fbox{2} \otimes \fbox{2} \in B^{\otimes 4}_1$, etc.
In the $i$-signature, one eliminates 
the neighboring pair $+-$ (not $-+$) successively to 
finally reach the pattern 
$\overbrace{-\cdots -}^\alpha\overbrace{+ \cdots +}^\beta$
called reduced $i$-signature.
The result is independent of the order of 
the eliminations when it can be done 
simultaneously in more than one places.
The reduced $i$-signature tells that 
$\varepsilon_i(b_i \otimes \cdots \otimes b_k)=\alpha$ and 
$\varphi_i(b_i \otimes \cdots \otimes b_k)=\beta$. 
In the above example, we get
\begin{equation*}
\underset{1-{\rm signature}}{{\underset{0-{\rm signature}}
{\phantom{1-signature}}}}
\underset{-+\phantom{+++}}{\underset{----\phantom{+}}{11112}} 
\otimes 
\underset{}{\underset{}{12}} \otimes 
\underset{}{\underset{}{2}} \otimes 
\underset{\phantom{--}++}{\underset{\phantom{--}++}{1122}}
\end{equation*}
Thus $\varepsilon_0 = 4,\, \varphi_0 = 2,\; 
\varepsilon_1=1$ and $\varphi_1 = 3$.
Finally 
${\tilde f}_i$ hits the component 
that is responsible for the leftmost $+$ 
in the reduced $i$-signature making it $-$.  
Similarly, ${\tilde e}_i$ hits the component 
corresponding to the rightmost $-$ 
in the reduced $i$-signature making it $+$.
If there is no such $+$ or $-$ to hit, the result 
of the action is $0$. 
The Weyl group operator 
$s_i$ acts so as to change the reduced $i$-signature  
$\overbrace{-\cdots -}^\alpha\overbrace{+ \cdots +}^\beta$
into 
$\overbrace{-\cdots -}^\beta\overbrace{+ \cdots +}^\alpha$.
In the above example, we have
\begin{align*}
p &= 
11112 \otimes 12 \otimes 2 \otimes 1122\\
{\tilde f}_0(p) &= 
11112 \otimes 12 \otimes 2 \otimes 1112\\
{\tilde f}_1(p) &= 
11122 \otimes 12 \otimes 2 \otimes 1122\\
{\tilde e}_0(p) &= 
11122 \otimes 12 \otimes 2 \otimes 1122\\
{\tilde e}_1(p) &= 
11111 \otimes 12 \otimes 2 \otimes 1122\\
s_0(p) &= 
11222 \otimes 12 \otimes 2 \otimes 1122\\
s_1(p) &= 
11122 \otimes 12 \otimes 2 \otimes 1222.
\end{align*}
For both $i=0$ and $1$, 
note that $\wt(s_i(p))=-\wt(p)$ for any $p$, 
and $s_i(p)=p$ if $\wt(p)=0$.
In order that ${\tilde e}_1p=0$ to hold for any  
$p \in B_{l_1}\ot \cdots \ot B_{l_k}$,
it is necessary and sufficient that 
the reduced $1$-signature of $p$ to become $+ \cdots +$.

\subsection{Combinatorial $R$}\label{subsec:r}
The crystal $B_l$ admits the affinization $\Aff(B_l)$.
It is the infinite set 
$\Aff(B_l) = \{ b[d] \mid b \in B_l, d \in \Z\}$
endowed with the crystal structure
${\tilde e}_i(b[d]) = (\et{i}b)[d+\delta_{i,0}]$,
${\tilde f}_i(b[d]) = (\ft{i}b)[d-\delta_{i,0}]$.
The element $b[d]$ will also be denoted by 
$\zeta^d b$ to save the space,
where $\zeta$ is an indeterminate.

The isomorphism of the affine crystal 
$\Aff(B_l) \otimes \Aff(B_k) \overset{\sim}{\rightarrow}
\Aff(B_k) \otimes \Aff(B_l)$ is the unique bijection that  
commutes with Kashiwara operators 
(up to a constant shift of $H$ below).
It is the $q=0$ analogue of the quantum $R$ and
called the combinatorial $R$.
Explicitly it is given by \cite{HHIKTT,Y}
\begin{equation*}
R(x[d] \otimes y[e]) 
= {\tilde y}[e-H(x\ot y)] \ot {\tilde x}[d+H(x\ot y)],
\end{equation*} 
where ${\tilde x}=({\tilde x}_i), {\tilde y} = ({\tilde y}_i)$
are given by
\begin{equation*}
\begin{split}
&{\tilde x}_i = x_i+Q_i(x,y)-Q_{i-1}(x,y),\quad 
{\tilde y}_i = y_i+Q_{i-1}(x,y)-Q_i(x,y),\\
&Q_i(x,y) = \min(x_{i+1}, y_i), \quad H(x\ot y) = Q_0(x,y). 
\end{split}
\end{equation*}
Here $x \otimes y \simeq {\tilde y}\otimes {\tilde x}$ 
under the isomorphism 
$B_l \otimes B_k \simeq B_k \otimes B_l$, 
which is also called (classical) combinatorial $R$.
$H$ is called the local energy function.
It is characterized by the recursion relation:
\begin{equation}\label{eq:hrec}
H(\et{i}(x\ot y))=\left\{%
\begin{array}{ll}
H(x\ot y)-1&\mbox{ if }i=0,\ \et{0}(x\ot y)=\et{0}x \ot y,\ 
\et{0}({\tilde y}\ot {\tilde x})=\et{0}{\tilde y} \ot {\tilde x},\\
H(x\ot y)+1&\mbox{ if }i=0,\ \et{0}(x\ot y)=x \ot \et{0}y,\ 
\et{0}({\tilde y}\ot {\tilde x})={\tilde y} \ot \et{0}{\tilde x},\\
H(x\ot y)&\mbox{ otherwise}.
\end{array}\right.
\end{equation}
together with the connectedness of
the crystal $B_l \ot B_k$.
(The original $H$ \cite{KMN} is $-H$ here.) 
$Q_0(x, y)$ is 
a solution of (\ref{eq:hrec})
normalized so as to attain the minimum at 
$Q_0(u_l\otimes u_k) = 0$ and ranges over 
$0 \le Q_0 \le \min(l,k)$ on $B_l \otimes B_k$.
Here,
\begin{equation}\label{eq:ul}
u_l = (l,0) = \fbox{1...1} \in B_l
\end{equation}
denotes the highest element
with respect to the $sl_2$ subcrystal concerning $\et{1}, \ft{1}$.
The invariance $Q_i(x \ot y) = Q_i({\tilde y} \ot {\tilde x})$ holds.
When $l=k$, the classical part of the combinatorial $R$ is trivial:
\begin{equation}\label{eq:tri}
R(\zeta^d x \otimes \zeta^e y) 
= \zeta^{e-H(x\ot y)}x \ot \zeta^{d+H(x\ot y)}y
\quad \hbox{on } B_l \ot B_l.
\end{equation}
The combinatorial $R$ has the following properties:
\begin{align}
&(\omega \otimes \omega)R = R(\omega \otimes \omega)
\quad \hbox{on } B_l \otimes B_k,\label{eq:dynkin}\\
&R \,\varrho = \varrho\,R \quad\hbox{on } B_l \otimes B_k,
\label{eq:rev}\\
&(1 \otimes R)(R \otimes 1)(1 \otimes R) =
(R \otimes 1)(1 \otimes R)(R \otimes 1)\quad
\hbox{on} \;
\Aff(B_j) \otimes \Aff(B_l) \otimes \Aff(B_k).\label{eq:ybr}
\end{align}
Here $\omega$ is the involutive automorphism (\ref{eq:odef}).
$\varrho(b_1 \otimes \cdots \otimes b_k) 
= b_k \otimes \cdots \otimes b_1$ is the reverse ordering 
of the tensor product for any $k$.
(\ref{eq:ybr}) is the Yang-Baxter relation.

To calculate the combinatorial $R$,
it is convenient to use the graphical rule 
(\cite{NY} Rule 3.11).
Consider the two elements $x=(x_1,x_2)\in B_k$
and $y=(y_1,y_2)\in B_l$.
Draw the following diagram to express
the tensor product $x\otimes y$.
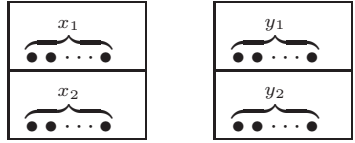
\begin{center}
\unitlength 13pt
\begin{picture}(10,4.5)
\multiput(0,0)(6,0){2}{
\multiput(0,0)(4,0){2}{\line(0,1){4}}
\multiput(0,0)(0,2){3}{\line(1,0){4}}
}
\put(0.5,0.2){$\overbrace{\bullet\bullet\cdots\bullet}^{x_{2}}$}
\put(0.5,2.2){$\overbrace{\bullet\bullet\cdots\bullet}^{x_1}$}
\put(6.5,0.2){$\overbrace{\bullet\bullet\cdots\bullet}^{y_{2}}$}
\put(6.5,2.2){$\overbrace{\bullet\bullet\cdots\bullet}^{y_1}$}
\end{picture}
\end{center}
The combinatorial $R$ and energy function $H$ for
$x\otimes y\in B_k\otimes B_l$ (with $k\geq l$) are found by
the following rule.
\begin{enumerate}
\item
Pick any dot, say $\bullet_a$, in the right column and connect it
with a dot $\bullet_a'$ in the left column by a line.
The partner $\bullet_a'$ is chosen
{}from the dots
whose positions are higher than that of $\bullet_a$.
If there is no such a dot, go to the bottom, and
the partner $\bullet_a'$ is chosen from the dots
in the lower row.
In the former case, we call such a pair ``unwinding,"
and, in the latter case, we call it ``winding."

\item
Repeat procedure (i) for the remaining unconnected dots
$(l-1)$ times.

\item
The isomorphism is obtained by
moving all unpaired dots in the left column to the right
horizontally.
We do not touch the paired dots during this move.

\item
The energy function $H$ is given by the number of unwinding pairs.
\end{enumerate}

The number of winding (unwinding) pairs is called
the winding (unwinding) number.
It is known that the result is independent of the
order of making pairs
(\cite{NY}, Propositions 3.15 and 3.17).
In the above description, we only consider the case $k\geq l$.
The other case $k\leq l$ can be done by reversing the above
procedure, noticing the fact $R^2=\mathrm{id}$.
For more properties, including that the above
definition indeed satisfies the axiom, see \cite{NY}.

\begin{example}
Corresponding to the tensor product
$\fbox{12222}\otimes\fbox{1122}$,
we draw diagram of the left hand side of
the following.
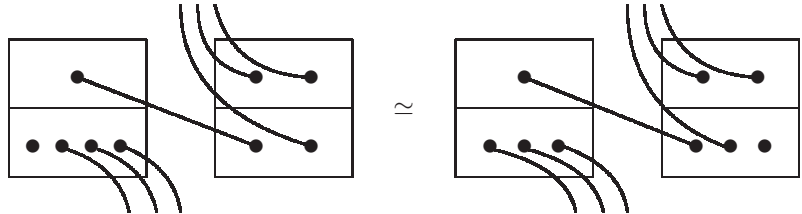
\begin{center}
\unitlength 13pt
\begin{picture}(24,6.5)(-0.2,-1)
\multiput(0,0)(6,0){2}{
\multiput(0,0)(4,0){2}{\line(0,1){4}}
\multiput(0,0)(0,2){3}{\line(1,0){4}}
}
\multiput(0.7,0.9)(0.85,0){4}{\circle*{0.4}}
\put(2.0,2.9){\circle*{0.4}}
\multiput(7.2,0.9)(1.6,0){2}{\circle*{0.4}}
\multiput(7.2,2.9)(1.6,0){2}{\circle*{0.4}}
\put(11.2,1.8){$\simeq$}
\thicklines
\qbezier(2.0,2.9)(2.0,2.9)(7.2,0.9)
\qbezier(8.8,0.9)(5,2)(5,5)
\qbezier(7.2,2.9)(5.5,3.2)(5.5,5)
\qbezier(8.8,2.9)(6.5,3)(6,5)
\qbezier(3.3,0.9)(4.8,0.3)(5,-1)
\qbezier(2.6,0.8)(4,0.3)(4.3,-1)
\qbezier(1.7,0.8)(3.2,0.3)(3.5,-1)
\thinlines
\put(13,0){
\multiput(0,0)(6,0){2}{
\multiput(0,0)(4,0){2}{\line(0,1){4}}
\multiput(0,0)(0,2){3}{\line(1,0){4}}
}
\multiput(1.0,0.9)(1.0,0){3}{\circle*{0.4}}
\put(2.0,2.9){\circle*{0.4}}
\multiput(7.0,0.9)(1.0,0){3}{\circle*{0.4}}
\multiput(7.2,2.9)(1.6,0){2}{\circle*{0.4}}
\thicklines
\qbezier(2.0,2.9)(2.0,2.9)(7.0,0.9)
\qbezier(7.8,0.9)(5,2)(5,5)
\qbezier(7.2,2.9)(5.5,3.2)(5.5,5)
\qbezier(8.8,2.9)(6.5,3)(6,5)
\qbezier(3.1,0.9)(4.8,0.3)(5,-1)
\qbezier(2.1,0.8)(4,0.3)(4.3,-1)
\qbezier(1.1,0.8)(3.2,0.3)(3.5,-1)
}
\end{picture}
\end{center}
By moving one unpaired dot to the right,
we obtain
\begin{equation*}
\fbox{12222}\otimes\fbox{1122}
\simeq
\fbox{1222}\otimes\fbox{11222}\,.
\end{equation*}
Since we have one unwinding pair, the energy function is
$H(\fbox{12222}\otimes\fbox{1122})=1$.

\end{example}

\section{Kerov--Kirillov--Reshetikhin bijection}\label{sec:kkr}
\subsection{Rigged configurations}
The Kerov-Kirillov-Reshetikhin (KKR) bijection 
gives a one to one correspondence between 
rigged configurations and highest paths.
The latter are elements of 
$B_{\lambda_1} \otimes \cdots \otimes B_{\lambda_L}$ 
annihilated by $\et{1}$. See around (\ref{eq:pdef}).

Let us define the rigged configurations.
Consider a pair 
$\left(\lambda ,\mu\right)$,
where both $\lambda$ and $\mu$ are positive integer sequences:
\begin{eqnarray*}
\lambda &=&
\left( \lambda_1,\lambda_2,\cdots ,\lambda_{L}\right),
\qquad (
L\in\mathbb{Z}_{\geq 0},\,
\lambda_i\in \mathbb{Z}_{>0}),\\
\mu &=&
\left( \mu_1,\mu_2,\cdots ,\mu_{N}\right),
\qquad (
N\in\mathbb{Z}_{\geq 0},\,
\mu_i\in \mathbb{Z}_{>0}),
\end{eqnarray*}
We use usual Young diagrammatic expression for
these integer sequences,
although $\lambda$, $\mu$ are not necessarily
assumed to be weakly decreasing.

\begin{definition}
$(1)$ For a given diagram $\nu$, we introduce coordinates
(row, column) of each boxes just like matrix entries.
For a box $\alpha$ of $\nu$, $\mathrm{col}(\alpha)$ is
the column coordinate of $\alpha$.
Then we define the following subsets:
\begin{equation*}
\nu |_{\leq j}:=\{ \alpha|\alpha\in\nu\,
,\mathrm{col}(\alpha )\leq j\},\quad
\nu |_{\geq j}:=\{ \alpha|\alpha\in\nu\,
,\mathrm{col}(\alpha )\geq j\}\, .
\end{equation*}

$(2)$ For the diagrams
$(\lambda,\mu)$,
we define $Q_j^{(a)}$ ($a=0,1$) by 
\begin{equation}\label{eq:q01}
Q_j^{(0)}:=\sum_{k=1}^{L}\min (j,\lambda_k),
\qquad
Q_j^{(1)}:=\sum_{k=1}^{N}\min (j,\mu_k),
\end{equation}
i.e., the number of boxes in $\lambda |_{\leq j}$
and $\mu |_{\leq j}$.
Then the vacancy number $p_j$ for rows of
$\mu$ is defined by 
\begin{equation}\label{eq:vn}
p_j:=Q^{(0)}_j-2Q^{(1)}_j,
\end{equation}
where $j$ is the width of the corresponding row.

\end{definition}

\begin{definition}\label{s:def:rc}
Consider the following data:
\begin{equation*}
\mathrm{RC}:=\bigl( \lambda,\,(\mu,r)\bigl) 
=\left( (\lambda_i)_{i=1}^L,\,
(\mu_i,r_i)_{i=1}^N\right) .
\end{equation*}

$(1)$ Calculate the vacancy numbers with respect
to the pair $(\lambda ,\mu)$.
If all vacancy numbers for rows of $\mu$
are nonnegative, i.e., 
$0\leq p_{\mu_i},
(1\leq i\leq N)$,
then RC is called a {\em configuration}.

$(2)$ If the integer $r_i$ satisfies the condition
\begin{equation}\label{eq:0rp}
0\leq r_i\leq p_{\mu_i},
\end{equation}
then $r_i$ is called a {\em rigging} associated with the row
$\mu_i$.
For the rows of equal widths, i.e., $\mu_i=\mu_{i+1}$,
we assume that $r_i\leq r_{i+1}$.

$(3)$ If RC is a configuration and if all integers $r_i$
are riggings associated with row $\mu_i$,
then RC is called $\mathfrak{sl}_2$ rigged configuration.

\end{definition}
In the rigged configuration, $\lambda$ is sometimes called a
quantum space which determines the shape of the corresponding path,
as we will see in the next subsection.
In the diagrammatic expression of rigged configurations,
the riggings are attached to the right of the
corresponding row.

\begin{definition}
For a given rigged configuration,
consider a row $\mu_i$ and the corresponding rigging $r_i$.
If they satisfy the condition
$r_i=p_{\mu_i}$,
then the row $\mu_i$ is called {\em singular}.

\end{definition}

\subsection{Definition of the KKR bijection}
Here we explain the original combinatorial
procedure to obtain a rigged configuration RC
\begin{equation*}
\phi :b
\longrightarrow {\rm RC}=\bigl( \lambda,\,(\mu,r)\bigl)
\end{equation*}
{}from a given path 
$b=b_1\otimes b_2\otimes\cdots\otimes b_L\in
B_{\lambda_1}\otimes B_{\lambda_2}\otimes\cdots\otimes B_{\lambda_L}
$.
The appearance of $\lambda_i$ in the right hand side
is clear from the following definition.

\begin{definition}\label{s:def:kkr}
Under the above setting, the image of the KKR map $\phi$ is defined
by the following procedure.

\begin{enumerate}
\item
We start from the empty rigged configuration
$\mathrm{RC}_0:=\bigl(\emptyset ,\, (\emptyset ,\emptyset )\bigl)$
and construct RC${}_1$, $\cdots$, RC${}_{|\lambda |}$
successively as follows (note that $|\lambda |=\sum\lambda_i$).

\item
Set 
$b_{1,0}:=b_1\in B_{\lambda_1}$
for the path $b=b_1\otimes \cdots$.
{}From $b_{1,0}$, we recursively construct
$b_{1,1}$, $b_{1,2}$, $\cdots$, $b_{1,\lambda_1}$
and RC${}_1$, RC${}_2$, $\cdots$, RC${}_{\lambda_1}$.
$b_{1,i+1}$ and RC${}_{i+1}$ are constructed from
$b_{1,i}:=(x_{1},x_{2})\in B_{\lambda_1-i}$
and RC${}_{i}$ as follows:

\begin{enumerate}
\item
First, assume that $x_{2}=0$.
Then we set
$b_{1,i+1}=(x_{1}-1,0)$.
If $i=0$, we create a new row to the quantum space as follows:
\begin{equation*}
\mathrm{RC}_1=\bigl( 1,\, (\emptyset ,\emptyset )\bigl) .
\end{equation*}
If $i>0$, then we add one box to the row of the quantum space
which is lengthened when we construct RC${}_{i}$.

\item
On the contrary, assume that $x_{2}>0$.
Then we set
$b_{1,i+1}=(x_{1},x_2-1)$.
We add a box to the quantum space by the same procedure
as in the case $x_{2}=0$.
Operation on $(\mu ,r)$ part of RC${}_{i}$ is as follows.
Calculate the vacancy numbers of RC${}_{i}$ and determine
all the singular rows.
If there is no singular row in $\mu |_{\geq i}$,
then create a new row in $\mu$.
On the other hand, assume that there is at least one singular row
in $\mu |_{\geq i}$.
Then, among these singular rows, we choose one of the longest
singular rows arbitrary, and let us tentatively call it $\mu_s$.
We add one box to the row $\mu_s$ and do not change the other parts,
which gives $\mu$ of new RC${}_{i+1}$.
As for the riggings, calculate the vacancy numbers of RC${}_{i+1}$.
Then we choose $r_s$, i.e., the rigging associated to
the lengthened row $\mu_s$,
so as to make the row $\mu_s$ singular in RC${}_{i+1}$.
Other riggings are chosen to be the same as
the corresponding ones in RC${}_{i}$.

\item
Repeat the above Step (b) for all letters 2
contained in $b_1$, then we do Step (a) for the rest of
letters 1 in $b_1$.
\end{enumerate}

\item
Do the same procedure of Step (ii) for $b_2$, $\cdots$, $b_L$.
Each time when we start with a new element $b_i$,
we create a new row in the quantum space,
and apply Step (ii).
The result gives  
the image RC=RC${}_{|\lambda |}$ of the map $\phi$.

\end{enumerate}
\end{definition}
It is known that all RC${}_i$ in the above procedure
are rigged configurations.

\begin{example}\label{s:ex:kkr}
Consider the following path:
\begin{equation*}
b=
\fbox{1111}\otimes\fbox{11}\otimes\fbox{22}\otimes
\fbox{12}\otimes\fbox{2}\otimes\fbox{122}\otimes
\fbox{122}\otimes\fbox{1112}
\end{equation*}
{}From the above path, we obtain the sequence of letters
$1111\cdot 11\cdot 22\cdot 21\cdot 2\cdot 221\cdot
221\cdot 2111$.
Then the calculation of $\phi (b)$ proceeds as follows.
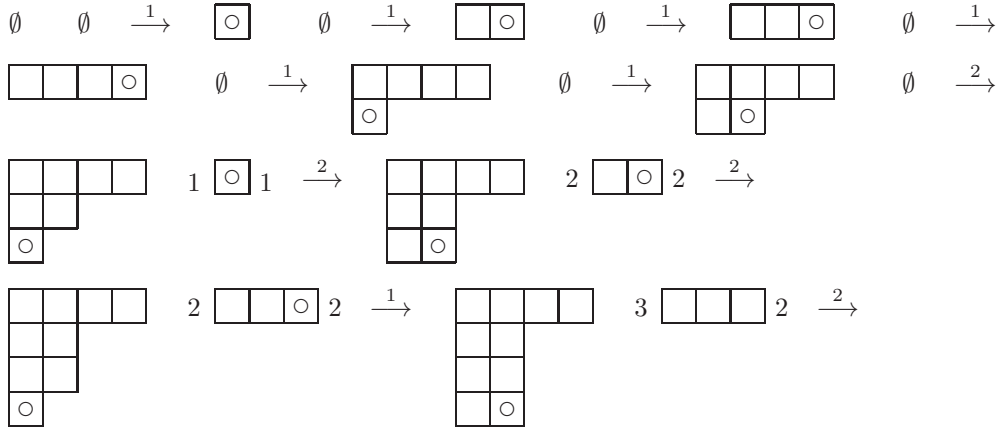
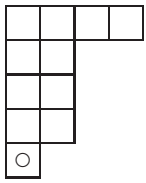
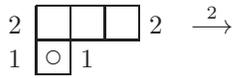
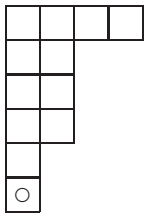
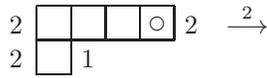
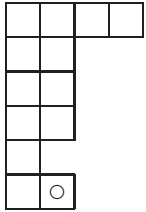
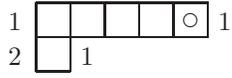
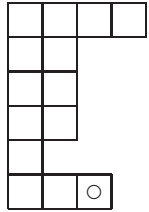
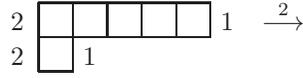
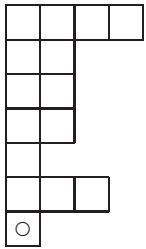
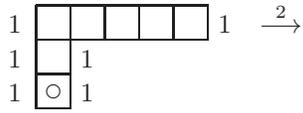
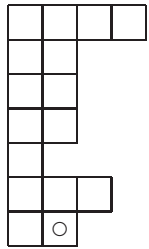
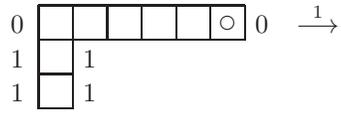
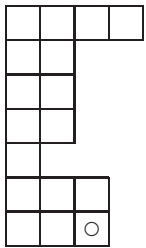
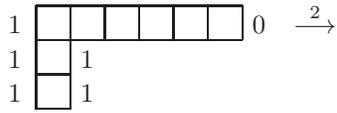
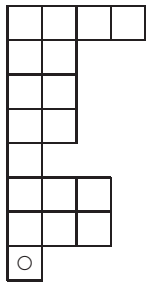
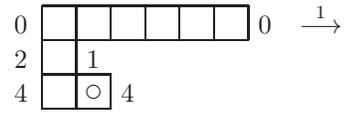
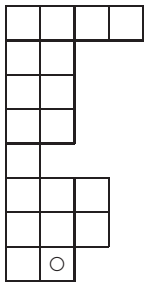
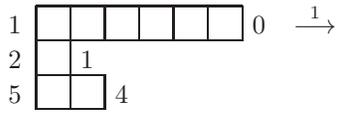
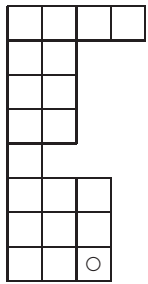
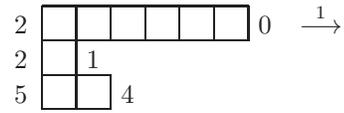
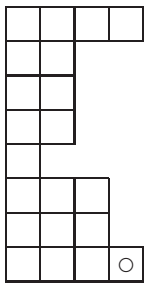
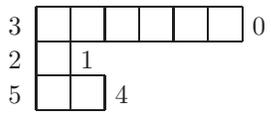
\begin{center}
\unitlength 13pt
\begin{picture}(30,1.5)(0,0)
\put(0,0.2){$\emptyset$}
\put(2,0.2){$\emptyset$}
\put(3.5,0.2){$\stackrel{1}{\longrightarrow}$}
\multiput(6,0)(1,0){2}{\line(0,1){1}}
\multiput(6,0)(0,1){2}{\line(1,0){1}}
\put(9,0.2){$\emptyset$}
\put(10.5,0.2){$\stackrel{1}{\longrightarrow}$}
\put(6.5,0.5){\circle{0.4}}
\multiput(13,0)(1,0){3}{\line(0,1){1}}
\multiput(13,0)(0,1){2}{\line(1,0){2}}
\put(17,0.2){$\emptyset$}
\put(18.5,0.2){$\stackrel{1}{\longrightarrow}$}
\put(14.5,0.5){\circle{0.4}}
\multiput(21,0)(1,0){4}{\line(0,1){1}}
\multiput(21,0)(0,1){2}{\line(1,0){3}}
\put(26,0.2){$\emptyset$}
\put(27.5,0.2){$\stackrel{1}{\longrightarrow}$}
\put(23.5,0.5){\circle{0.4}}
\end{picture}
\end{center}
\begin{center}
\unitlength 13pt
\begin{picture}(30,2.5)(0,0)
\multiput(0,1)(1,0){5}{\line(0,1){1}}
\multiput(0,1)(0,1){2}{\line(1,0){4}}
\put(6,1.2){$\emptyset$}
\put(7.5,1.2){$\stackrel{1}{\longrightarrow}$}
\put(3.5,1.5){\circle{0.4}}
\multiput(10,1)(1,0){5}{\line(0,1){1}}
\multiput(10,1)(0,1){2}{\line(1,0){4}}
\multiput(10,0)(1,0){2}{\line(0,1){1}}
\put(10,0){\line(1,0){1}}
\put(16,1.2){$\emptyset$}
\put(17.5,1.2){$\stackrel{1}{\longrightarrow}$}
\put(10.5,0.5){\circle{0.4}}
\multiput(20,1)(1,0){5}{\line(0,1){1}}
\multiput(20,1)(0,1){2}{\line(1,0){4}}
\multiput(20,0)(1,0){3}{\line(0,1){1}}
\put(20,0){\line(1,0){2}}
\put(26,1.2){$\emptyset$}
\put(27.5,1.2){$\stackrel{2}{\longrightarrow}$}
\put(21.5,0.5){\circle{0.4}}
\end{picture}
\end{center}
\begin{center}
\unitlength 13pt
\begin{picture}(30,4.5)(0,-1)
\put(0,0){\line(1,0){1}}
\put(0,1){\line(1,0){2}}
\multiput(0,0)(1,0){2}{\line(0,1){1}}
\multiput(0,1)(1,0){3}{\line(0,1){1}}
\multiput(0,2)(1,0){5}{\line(0,1){1}}
\multiput(0,2)(0,1){2}{\line(1,0){4}}
\multiput(6,2)(1,0){2}{\line(0,1){1}}
\multiput(6,2)(0,1){2}{\line(1,0){1}}
\put(5.2,2.1){1}
\put(7.3,2.1){1}
\put(0.5,0.5){\circle{0.4}}
\put(6.5,2.5){\circle{0.4}}
\put(8.5,2.2){$\stackrel{2}{\longrightarrow}$}
\multiput(11,0)(1,0){3}{\line(0,1){2}}
\multiput(11,0)(0,1){2}{\line(1,0){2}}
\multiput(11,2)(1,0){5}{\line(0,1){1}}
\multiput(11,2)(0,1){2}{\line(1,0){4}}
\put(12.5,0.5){\circle{0.4}}
\multiput(17,2)(1,0){3}{\line(0,1){1}}
\multiput(17,2)(0,1){2}{\line(1,0){2}}
\put(18.5,2.5){\circle{0.4}}
\put(16.2,2.2){2}
\put(19.3,2.2){2}
\put(20.5,2.2){$\stackrel{2}{\longrightarrow}$}
\end{picture}
\end{center}
\begin{center}
\unitlength 13pt
\begin{picture}(30,4.5)(0,-1)
\put(-23,1){
\put(23,-1){\line(1,0){1}}
\multiput(23,-1)(1,0){2}{\line(0,1){1}}
\multiput(23,0)(1,0){3}{\line(0,1){2}}
\multiput(23,0)(0,1){2}{\line(1,0){2}}
\multiput(23,2)(1,0){5}{\line(0,1){1}}
\multiput(23,2)(0,1){2}{\line(1,0){4}}
\multiput(29,2)(1,0){4}{\line(0,1){1}}
\multiput(29,2)(0,1){2}{\line(1,0){3}}
\put(28.2,2.2){2}
\put(32.3,2.2){2}
\put(23.5,-0.5){\circle{0.4}}
\put(31.5,2.5){\circle{0.4}}
\put(33.5,2.2){$\stackrel{1}{\longrightarrow}$}
}
\put(13,-2){
\multiput(0,2)(1,0){3}{\line(0,1){3}}
\multiput(0,2)(0,1){3}{\line(1,0){2}}
\multiput(0,5)(1,0){5}{\line(0,1){1}}
\multiput(0,5)(0,1){2}{\line(1,0){4}}
\put(1.5,2.5){\circle{0.4}}
\multiput(6,5)(1,0){4}{\line(0,1){1}}
\multiput(6,5)(0,1){2}{\line(1,0){3}}
}
\put(18.2,3.2){3}
\put(22.3,3.2){2}
\put(23.5,3.2){$\stackrel{2}{\longrightarrow}$}
\end{picture}
\end{center}
\begin{center}
\unitlength 13pt
\begin{picture}(30,5.5)
\put(0,1){\line(1,0){1}}
\multiput(0,1)(1,0){2}{\line(0,1){1}}
\multiput(0,2)(1,0){3}{\line(0,1){3}}
\multiput(0,2)(0,1){3}{\line(1,0){2}}
\multiput(0,5)(1,0){5}{\line(0,1){1}}
\multiput(0,5)(0,1){2}{\line(1,0){4}}
\put(0.5,1.5){\circle{0.4}}
\put(6,4){\line(1,0){1}}
\multiput(6,4)(1,0){2}{\line(0,1){1}}
\multiput(6,5)(1,0){4}{\line(0,1){1}}
\multiput(6,5)(0,1){2}{\line(1,0){3}}
\put(6.5,4.5){\circle{0.4}}
\put(5.2,4.2){1}
\put(5.2,5.2){2}
\put(7.3,4.2){1}
\put(9.3,5.2){2}
\put(10.5,5.2){$\stackrel{2}{\longrightarrow}$}
\multiput(13,0)(1,0){2}{\line(0,1){2}}
\multiput(13,0)(0,1){2}{\line(1,0){1}}
\multiput(13,2)(1,0){3}{\line(0,1){3}}
\multiput(13,2)(0,1){3}{\line(1,0){2}}
\multiput(13,5)(1,0){5}{\line(0,1){1}}
\multiput(13,5)(0,1){2}{\line(1,0){4}}
\put(13.5,0.5){\circle{0.4}}
\put(19,4){\line(1,0){1}}
\multiput(19,4)(1,0){2}{\line(0,1){1}}
\multiput(19,5)(1,0){5}{\line(0,1){1}}
\multiput(19,5)(0,1){2}{\line(1,0){4}}
\put(22.5,5.5){\circle{0.4}}
\put(18.2,4.2){2}
\put(18.2,5.2){2}
\put(20.3,4.2){1}
\put(23.3,5.2){2}
\put(24.5,5.2){$\stackrel{2}{\longrightarrow}$}
\end{picture}
\end{center}
\begin{center}
\unitlength 13pt
\begin{picture}(30,6.5)
\multiput(0,0)(1,0){2}{\line(0,1){2}}
\multiput(0,0)(0,1){2}{\line(1,0){2}}
\put(2,0){\line(0,1){1}}
\multiput(0,2)(1,0){3}{\line(0,1){3}}
\multiput(0,2)(0,1){3}{\line(1,0){2}}
\multiput(0,5)(1,0){5}{\line(0,1){1}}
\multiput(0,5)(0,1){2}{\line(1,0){4}}
\put(1.5,0.5){\circle{0.4}}
\put(6,4){\line(1,0){1}}
\multiput(6,4)(1,0){2}{\line(0,1){1}}
\multiput(6,5)(1,0){6}{\line(0,1){1}}
\multiput(6,5)(0,1){2}{\line(1,0){5}}
\put(10.5,5.5){\circle{0.4}}
\put(5.2,4.2){2}
\put(5.2,5.2){1}
\put(7.3,4.2){1}
\put(11.3,5.2){1}
\put(12.5,5.2){$\stackrel{1}{\longrightarrow}$}
\multiput(15,0)(1,0){2}{\line(0,1){2}}
\multiput(15,0)(0,1){2}{\line(1,0){3}}
\multiput(17,0)(1,0){2}{\line(0,1){1}}
\multiput(15,2)(1,0){3}{\line(0,1){3}}
\multiput(15,2)(0,1){3}{\line(1,0){2}}
\multiput(15,5)(1,0){5}{\line(0,1){1}}
\multiput(15,5)(0,1){2}{\line(1,0){4}}
\put(17.5,0.5){\circle{0.4}}
\put(21,4){\line(1,0){1}}
\multiput(21,4)(1,0){2}{\line(0,1){1}}
\multiput(21,5)(1,0){6}{\line(0,1){1}}
\multiput(21,5)(0,1){2}{\line(1,0){5}}
\put(20.2,4.2){2}
\put(20.2,5.2){2}
\put(22.3,4.2){1}
\put(26.3,5.2){1}
\put(27.5,5.2){$\stackrel{2}{\longrightarrow}$}
\end{picture}
\end{center}
\begin{center}
\unitlength 13pt
\begin{picture}(30,7.5)
\multiput(0,1)(1,0){2}{\line(0,1){2}}
\multiput(0,1)(0,1){2}{\line(1,0){3}}
\multiput(2,1)(1,0){2}{\line(0,1){1}}
\multiput(0,3)(1,0){3}{\line(0,1){3}}
\multiput(0,3)(0,1){3}{\line(1,0){2}}
\multiput(0,6)(1,0){5}{\line(0,1){1}}
\multiput(0,6)(0,1){2}{\line(1,0){4}}
\put(0,0){\line(1,0){1}}
\multiput(0,0)(1,0){2}{\line(0,1){1}}
\put(0.5,0.5){\circle{0.4}}
\multiput(6,4)(0,1){2}{\line(1,0){1}}
\multiput(6,4)(1,0){2}{\line(0,1){2}}
\multiput(6,6)(1,0){6}{\line(0,1){1}}
\multiput(6,6)(0,1){2}{\line(1,0){5}}
\put(6.5,4.5){\circle{0.4}}
\put(5.2,4.2){1}
\put(5.2,5.2){1}
\put(5.2,6.2){1}
\put(7.3,4.2){1}
\put(7.3,5.2){1}
\put(11.3,6.2){1}
\put(12.5,6.2){$\stackrel{2}{\longrightarrow}$}
\multiput(15,1)(1,0){2}{\line(0,1){2}}
\multiput(15,1)(0,1){2}{\line(1,0){3}}
\multiput(17,1)(1,0){2}{\line(0,1){1}}
\multiput(15,3)(1,0){3}{\line(0,1){3}}
\multiput(15,3)(0,1){3}{\line(1,0){2}}
\multiput(15,6)(1,0){5}{\line(0,1){1}}
\multiput(15,6)(0,1){2}{\line(1,0){4}}
\put(15,0){\line(1,0){2}}
\multiput(15,0)(1,0){3}{\line(0,1){1}}
\put(16.5,0.5){\circle{0.4}}
\multiput(21,4)(0,1){2}{\line(1,0){1}}
\multiput(21,4)(1,0){2}{\line(0,1){2}}
\multiput(21,6)(1,0){7}{\line(0,1){1}}
\multiput(21,6)(0,1){2}{\line(1,0){6}}
\put(26.5,6.5){\circle{0.4}}
\put(20.2,4.2){1}
\put(20.2,5.2){1}
\put(20.2,6.2){0}
\put(22.3,4.2){1}
\put(22.3,5.2){1}
\put(27.3,6.2){0}
\put(28.5,6.2){$\stackrel{1}{\longrightarrow}$}
\end{picture}
\end{center}
\begin{center}
\unitlength 13pt
\begin{picture}(30,8.5)
\multiput(0,2)(1,0){2}{\line(0,1){2}}
\multiput(0,2)(0,1){2}{\line(1,0){3}}
\multiput(2,2)(1,0){2}{\line(0,1){1}}
\multiput(0,4)(1,0){3}{\line(0,1){3}}
\multiput(0,4)(0,1){3}{\line(1,0){2}}
\multiput(0,7)(1,0){5}{\line(0,1){1}}
\multiput(0,7)(0,1){2}{\line(1,0){4}}
\put(0,1){\line(1,0){3}}
\multiput(0,1)(1,0){4}{\line(0,1){1}}
\put(2.5,1.5){\circle{0.4}}
\multiput(6,5)(0,1){2}{\line(1,0){1}}
\multiput(6,5)(1,0){2}{\line(0,1){2}}
\multiput(6,7)(1,0){7}{\line(0,1){1}}
\multiput(6,7)(0,1){2}{\line(1,0){6}}
\put(5.2,5.2){1}
\put(5.2,6.2){1}
\put(5.2,7.2){1}
\put(7.3,5.2){1}
\put(7.3,6.2){1}
\put(12.3,7.2){0}
\put(13.5,7.2){$\stackrel{2}{\longrightarrow}$}
\put(16,0){
\multiput(0,2)(1,0){2}{\line(0,1){2}}
\multiput(0,2)(0,1){2}{\line(1,0){3}}
\multiput(2,2)(1,0){2}{\line(0,1){1}}
\multiput(0,4)(1,0){3}{\line(0,1){3}}
\multiput(0,4)(0,1){3}{\line(1,0){2}}
\multiput(0,7)(1,0){5}{\line(0,1){1}}
\multiput(0,7)(0,1){2}{\line(1,0){4}}
\put(0,1){\line(1,0){3}}
\multiput(0,1)(1,0){4}{\line(0,1){1}}
\put(0,0){\line(1,0){1}}
\multiput(0,0)(1,0){2}{\line(0,1){1}}
\put(0.5,0.5){\circle{0.4}}
\multiput(6,5)(0,1){2}{\line(1,0){1}}
\multiput(6,5)(1,0){2}{\line(0,1){2}}
\put(7,6){\line(1,0){1}}
\put(7,5){\line(1,0){1}}
\put(8,5){\line(0,1){1}}
\put(7.5,5.5){\circle{0.4}}
\multiput(6,7)(1,0){7}{\line(0,1){1}}
\multiput(6,7)(0,1){2}{\line(1,0){6}}
\put(5.2,5.2){4}
\put(5.2,6.2){2}
\put(5.2,7.2){0}
\put(8.3,5.2){4}
\put(7.3,6.2){1}
\put(12.3,7.2){0}
\put(13.5,7.2){$\stackrel{1}{\longrightarrow}$}
}
\end{picture}
\end{center}
\begin{center}
\unitlength 13pt
\begin{picture}(30,8.5)
\multiput(0,2)(1,0){2}{\line(0,1){2}}
\multiput(0,2)(0,1){2}{\line(1,0){3}}
\multiput(2,2)(1,0){2}{\line(0,1){1}}
\multiput(0,4)(1,0){3}{\line(0,1){3}}
\multiput(0,4)(0,1){3}{\line(1,0){2}}
\multiput(0,7)(1,0){5}{\line(0,1){1}}
\multiput(0,7)(0,1){2}{\line(1,0){4}}
\put(0,1){\line(1,0){3}}
\multiput(0,1)(1,0){4}{\line(0,1){1}}
\put(0,0){\line(1,0){2}}
\multiput(0,0)(1,0){3}{\line(0,1){1}}
\put(1.5,0.5){\circle{0.4}}
\multiput(6,5)(0,1){2}{\line(1,0){1}}
\multiput(6,5)(1,0){2}{\line(0,1){2}}
\put(7,6){\line(1,0){1}}
\put(7,5){\line(1,0){1}}
\put(8,5){\line(0,1){1}}
\multiput(6,7)(1,0){7}{\line(0,1){1}}
\multiput(6,7)(0,1){2}{\line(1,0){6}}
\put(5.2,5.2){5}
\put(5.2,6.2){2}
\put(5.2,7.2){1}
\put(8.3,5.2){4}
\put(7.3,6.2){1}
\put(12.3,7.2){0}
\put(13.5,7.2){$\stackrel{1}{\longrightarrow}$}
\put(16,0){
\multiput(0,2)(1,0){2}{\line(0,1){2}}
\multiput(0,2)(0,1){2}{\line(1,0){3}}
\multiput(2,2)(1,0){2}{\line(0,1){1}}
\multiput(0,4)(1,0){3}{\line(0,1){3}}
\multiput(0,4)(0,1){3}{\line(1,0){2}}
\multiput(0,7)(1,0){5}{\line(0,1){1}}
\multiput(0,7)(0,1){2}{\line(1,0){4}}
\put(0,1){\line(1,0){3}}
\multiput(0,1)(1,0){4}{\line(0,1){1}}
\put(0,0){\line(1,0){3}}
\multiput(0,0)(1,0){4}{\line(0,1){1}}
\put(2.5,0.5){\circle{0.4}}
\multiput(6,5)(0,1){2}{\line(1,0){1}}
\multiput(6,5)(1,0){2}{\line(0,1){2}}
\put(7,6){\line(1,0){1}}
\put(7,5){\line(1,0){1}}
\put(8,5){\line(0,1){1}}
\multiput(6,7)(1,0){7}{\line(0,1){1}}
\multiput(6,7)(0,1){2}{\line(1,0){6}}
\put(5.2,5.2){5}
\put(5.2,6.2){2}
\put(5.2,7.2){2}
\put(8.3,5.2){4}
\put(7.3,6.2){1}
\put(12.3,7.2){0}
\put(13.5,7.2){$\stackrel{1}{\longrightarrow}$}
}
\end{picture}
\end{center}
\begin{center}
\unitlength 13pt
\begin{picture}(30,8.5)
\multiput(0,2)(1,0){2}{\line(0,1){2}}
\multiput(0,2)(0,1){2}{\line(1,0){3}}
\multiput(2,2)(1,0){2}{\line(0,1){1}}
\multiput(0,4)(1,0){3}{\line(0,1){3}}
\multiput(0,4)(0,1){3}{\line(1,0){2}}
\multiput(0,7)(1,0){5}{\line(0,1){1}}
\multiput(0,7)(0,1){2}{\line(1,0){4}}
\put(0,1){\line(1,0){3}}
\multiput(0,1)(1,0){4}{\line(0,1){1}}
\put(0,0){\line(1,0){4}}
\multiput(0,0)(1,0){5}{\line(0,1){1}}
\put(3,1){\line(1,0){1}}
\put(3.5,0.5){\circle{0.4}}
\multiput(6,5)(0,1){2}{\line(1,0){1}}
\multiput(6,5)(1,0){2}{\line(0,1){2}}
\put(7,6){\line(1,0){1}}
\put(7,5){\line(1,0){1}}
\put(8,5){\line(0,1){1}}
\multiput(6,7)(1,0){7}{\line(0,1){1}}
\multiput(6,7)(0,1){2}{\line(1,0){6}}
\put(5.2,5.2){5}
\put(5.2,6.2){2}
\put(5.2,7.2){3}
\put(8.3,5.2){4}
\put(7.3,6.2){1}
\put(12.3,7.2){0}
\end{picture}
\end{center}
In the above diagrams, newly added boxes are
indicated by circles ``$\circ$",
and vacancy numbers are listed on the left
of the corresponding rows in order to
facilitate the calculations.
This example, fully worked out here, will be revisited in
Example \ref{s:ex:mainth} by using Theorem \ref{s:main}
for comparison.

\end{example}

\subsection{Basic properties of the KKR bijection}
It is known that the inverse map $\phi ^{-1}$
can be described by a similar combinatorial procedure.
We will use both $\phi$ and $\phi^{-1}$ in later arguments.

\begin{theorem}\label{s:def:rkk}
The inverse map
\begin{equation*}
\phi^{-1} :{\rm RC}=\bigl( \lambda,\,(\mu,r)\bigl)
\longrightarrow b,
\end{equation*}
is obtained by the following procedure
($\lambda =(\lambda_1,\lambda_2,\cdots ,\lambda_L)$).
\begin{enumerate}
\item
We construct RC${}_{|\lambda |}$,
RC${}_{|\lambda |-1}$, $\cdots$,
RC${}_{1}$,
RC${}_{0}=\bigl(\emptyset ,\, (\emptyset ,\emptyset )\bigl)$ 
and $b_L$, $b_{L-1}$, $\cdots$, $b_1$ ($b_i\in B_{\lambda_i}$)
as follows.
\item
We start from $\lambda_L$ of $\lambda$,
and set $\mathrm{RC}_{|\lambda |}:=\mathrm{RC}$
and $b_{L,\lambda_L}=(0,0)$.
In order to obtain the quantum space of
RC${}_{|\lambda |-1}$, RC${}_{|\lambda |-2}$, $\cdots$,
we remove boxes of row $\lambda_L$ of the quantum
space one by one.
RC${}_{i-1}$ and $b_{L,i-1}$ are constructed from
RC${}_{i}$ and $b_{L,i}$.
We call the rightmost box of the row of 
length $\lambda_L-i$ in the quantum space as $\alpha$.
(${\rm col}(\alpha) = \lambda_L-i$.)
Let us tentatively denote the ``$\mu$ part" of
RC${}_{i}$ by $\nu$.
\begin{enumerate}
\item
Assume that there is no singular row in
$\nu |_{\geq\mathrm{col}(\alpha )}$.
Then RC${}_{i-1}$ is obtained by  removing
the box $\alpha$ from the quantum space.
In this case, $b_{L,i-1}$ is obtained by
adding letter 1 to $b_{L,i}$.
\item
Assume on the contrary that there are singular rows
in $\nu |_{\geq\mathrm{col}(\alpha )}$.
Among these singular rows, we choose one of the
shortest rows arbitrary, and denote the
rightmost box of the chosen row by $\beta$.
Then RC${}_{i-1}$ is obtained by removing the
two boxes $\alpha$ and $\beta$ from RC${}_{i}$.
New riggings are specified as follows.
For the row from which the box $\beta$ is removed,
take new rigging so that it becomes
singular in RC${}_{i-1}$.
On the contrary, for the other rows 
riggings are kept unchanged from 
the corresponding ones in RC${}_{i}$.
Finally, $b_{L,i-1}$ is obtained by
adding letter 2 to $b_{L,i}$.
\item
By doing Steps (a) and (b) for the rest of
boxes of $\lambda_L$ in the quantum space,
we obtain $b_L\in B_{\lambda_L}$.
Here, the orderings of letters 1 and 2 within
$b_L$ is chosen so that $b_L$ becomes
semi-standard Young tableau.
\end{enumerate}
\item
By doing Step (ii) for the rest of rows
$\lambda_{L-1}$, $\cdots$, $\lambda_1$,
we obtain $b_{L-1}$, $\cdots$, $b_1$
respectively.
Finally, we obtain the path
$b=b_1\otimes b_2\otimes\cdots\otimes b_L$
as the image of the map $\phi^{-1}$.

\end{enumerate}
\end{theorem}

\begin{example}
For an example of the calculation of $\phi^{-1}$,
one can use Example \ref{s:ex:kkr},
that is, reverse all arrows ``$\longrightarrow$"
to ``$\longleftarrow$".

\end{example}

The above map $\phi^{-1}$ depends on 
ordering of the quantum space
$\lambda =(\lambda_1,\lambda_2,\cdots ,\lambda_L)$.
The dependence is described by 

\begin{theorem}[\cite{KSS} Lemma 8.5]\label{s:thm:KSS}
Let $\alpha, \beta$ be any successive two rows
in the quantum space of a rigged configuration.
Suppose the removal of $\alpha$ first
and $\beta$ next by $\phi^{-1}$
lead to the tableaux $a_1$ and $b_1$, respectively.
Suppose similarly that 
the removal of $\beta$ first
and $\alpha$ next lead to $b_2$ and $a_2$, respectively.
If the order of the other removal is the same, 
\begin{equation*}
b_1\otimes a_1\,\simeq\, a_2\otimes b_2
\end{equation*}
is valid under the isomorphism of the combinatorial $R$.
\end{theorem}

The KKR bijection $\phi^{\pm 1}$, originally 
designed only for highest paths, 
is known to admit an extension which covers
{\em all} the paths.
In fact one can apply the same combinatorial procedure
as $\phi$ to obtain $\phi (b)$ for any 
$b \in B_{\lambda_1}\otimes\cdots\otimes B_{\lambda_L}$.
The resulting object is a {\it unrestricted rigged configuration},
where the condition (\ref{eq:0rp}) is relaxed to 
$-\mu_i \le r_i \le p_{\mu_i}$ \cite{Sch2,DS}. 
Obviously rigged configurations are special case of 
unrestricted ones.
For unrestricted rigged configurations,
combinatorial procedure in Theorem \ref{s:def:rkk}
also works to define the inverse map $\phi^{-1}$.
Let us write $\phi(b)=(\lambda,(\mu,r))$.
Then, $|\lambda |$ represents the total number of letters 1 and
2 contained in the path $b$, whereas $|\mu |$
is the number of letter 2 in $b$.
Note in particular that $|\lambda |\geq |\mu |$
holds for unrestricted rigged configurations.

Given a non-highest path $b$, one can always make
\begin{equation*}
b':=\fbox{1}^{\,\otimes \Lambda}\otimes b
\end{equation*}
highest by taking 
$\Lambda\geq \lambda_1+\cdots +\lambda_L$.
Under these notations, we have 

\begin{lemma}\label{s:lem:haiesutoka}
Let the unrestricted rigged configuration
corresponding to $b$ be
\begin{equation}\label{s:eq:unristrictedrc}
\bigl( (\lambda_i)_{i=1}^L ,(\mu_j ,r_j)_{j=1}^N\bigl).
\end{equation}
Then the rigged configuration corresponding to
the highest path $b'$ is given by
\begin{equation}\label{s:eq:unristrictedrc_modify}
\bigl( (\lambda_i)_{i=1}^L\cup (1^\Lambda)
,(\mu_j ,r_j+\Lambda)_{j=1}^N\bigl).
\end{equation}
\end{lemma}
{\bf Proof.}
Let the vacancy number of a row $\mu_j$
of the pair $(\lambda ,\mu )$ of (\ref{s:eq:unristrictedrc})
be $p_{\mu_j}$.
Then the vacancy number of the row $\mu_j$ of
(\ref{s:eq:unristrictedrc_modify})
is equal to $p_{\mu_j}+\Lambda$.
Now we apply $\phi^{-1}$ on
(\ref{s:eq:unristrictedrc_modify}).
{}From the quantum space 
$\lambda\cup (1^\Lambda)$,
we remove $\lambda$ first, and next remove $(1^\Lambda)$.
Recall that the combinatorial procedure 
in Theorem \ref{s:def:rkk} only refers to co-rigging
(:=vacancy number $-$ rigging),
rather than the riggings.
Therefore, when we remove $\lambda$ from the quantum space
of (\ref{s:eq:unristrictedrc_modify}), we get
$b$ as the corresponding part of the image.
Then, the remaining rigged configuration has
the quantum space $(1^\Lambda)$ without $\mu$ part.
Since the map $\phi^{-1}$ becomes trivial on it, we obtain  
$b'$ as the image of (\ref{s:eq:unristrictedrc_modify}).

\begin{acknowledgments}
The authors thank Rei Inoue, 
Tomoki Nakanishi, Akira Takenouchi for discussion and 
Taichiro Takagi for collaboration in an early stage 
of the work and for communicating the result in \cite{Ta}.
A.K. is supported by Grants-in-Aid for Scientific No.19540393.
R.S. is a research fellow of the 
Japan Society for the Promotion of Science. 
\end{acknowledgments}


\vspace{5mm}
\begin{flushleft}
Atsuo Kuniba:\\
Institute of Physics, Graduate School of Arts and Sciences,
University of Tokyo,
Komaba, Tokyo 153-8902, Japan\\
\texttt{atsuo@gokutan.c.u-tokyo.ac.jp}\vspace{3mm}\\
Reiho Sakamoto:\\
Department of Physics, Graduate School of Science, 
University of Tokyo, Hongo, Tokyo 113-0033, Japan\\
\texttt{reiho@monet.phys.s.u-tokyo.ac.jp}\vspace{3mm}
\end{flushleft}

\end{document}